\documentclass[11pt]{article}
\usepackage[total={6in, 7.75in}]{geometry}
\usepackage{authblk}
\usepackage{amsmath,amssymb,amsthm,graphicx,tikz,multirow,hhline,color,array,dsfont}
\usepackage[labelfont=bf]{caption}
\usepackage[T1]{fontenc}
\usepackage{esint}

\usepackage[sort,numbers]{natbib}
\usepackage{xcolor}

\usepackage{hyperref}			
\hypersetup{colorlinks=true,linkcolor=black,urlcolor=blue,citecolor=blue}

\newtheorem{conj}{Conjecture}

\newtheorem{prop}{Proposition}
\theoremstyle{definition}
\newtheorem{remark}{Remark}

\newcommand{\beq}{\begin{equation}}
\newcommand{\eeq}{\end{equation}}

\newcommand{\ov}{\overline}
\newcommand{\ddt}{\frac{{\rm d}}{{\rm d}t}}
\newcommand{\dVdq}{\frac{\partial V}{\partial q^2}}
\newcommand{\dPdq}{\frac{\partial P}{\partial q^2}}
\newcommand{\ba}{\mathbf{a}}
\newcommand{\bb}{\mathbf{b}}
\newcommand{\f}{\mathbf{f}}
\newcommand{\bTh}{\mathbf{\Theta}}
\newcommand{\bM}{\mathbf{M}}
\newcommand{\E}{\mathcal{E}}
\newcommand{\LL}{W}
\newcommand{\U}{\mathcal{U}}
\newcommand{\V}{\mathcal{V}}

\newcommand{\R}{\mathbb{R}}
\newcommand{\dx}{{\rm d}x}

\newcommand\solidrule[1][15pt]{\rule[0.5ex]{#1}{1pt}}
\newcommand\dashedrule{\mbox{%
	\solidrule[3pt]\hspace{3pt}\solidrule[3pt]\hspace{3pt}\solidrule[3pt]}}
\newcommand\longdashedrule{\mbox{%
	\solidrule[6pt]\hspace{2pt}\solidrule[6pt]}}
\newcommand\dottedrule{\mbox{%
	\solidrule[1.5pt]\hspace{2pt}\solidrule[1.5pt]\hspace{2pt}\solidrule[1.5pt]\hspace{2pt}\solidrule[1.5pt]\hspace{2pt}\solidrule[1.5pt]}}

\title{Bounds on mean energy in the Kuramoto--Sivashinsky equation computed using semidefinite programming}
\author{David Goluskin$^{1,2,}$\thanks{Email: {\tt goluskin@uvic.ca}}~
	and Giovanni Fantuzzi$^3$}
\affil{$^1$Department of Mathematics and Statistics, University of Victoria, Victoria, BC, V8P 5C2, Canada\\
$^2$Department of Mathematics, University of Michigan, Ann Arbor, MI 48109, USA\\
$^3$Department of Aeronautics, Imperial College London, South Kensington Campus, London, SW7 2AZ, United Kingdom}
\date{}

\begin{document}

\maketitle

\begin{abstract}
We present methods for bounding infinite-time averages in dynamical systems governed by nonlinear PDEs. The methods rely on auxiliary functionals, which are similar to Lyapunov functionals but satisfy different inequalities. The inequalities are enforced by requiring certain expressions to be sums of squares of polynomials, and the optimal choice of auxiliary functional is posed as a semidefinite program (SDP) that can be solved computationally. To formulate these SDPs we approximate the PDE by truncated systems of ODEs and proceed in one of two ways. The first approach is to compute bounds for the ODE systems, increasing the truncation order until bounds converge numerically. The second approach incorporates the ODE systems with analytical estimates on their deviation from the PDE, thereby using finite truncations to produce bounds for the full PDE. We apply both methods to the Kuramoto--Sivashinsky equation. In particular, we compute upper bounds on the spatiotemporal average of energy by employing polynomial auxiliary functionals up to degree six. The first approach is used for most computations, but a subset of results are checked using the second approach, and the results agree to high precision. These bounds apply to all odd solutions of period $2\pi L$, where $L$ is varied. Sharp bounds are obtained for $L\le10$, and trends suggest that more expensive computations would yield sharp bounds at larger $L$ also. The bounds are known to be sharp (to within 0.1\% numerical error) because they are saturated by the simplest nonzero steady states, which apparently have the largest mean energy among all odd solutions. Prior authors have conjectured that mean energy remains $O(1)$ for $L\gg1$ since no particular solutions with larger energy have been found. Our bounds constitute the first positive evidence for this conjecture, albeit up to finite $L$, and they offer some guidance for analytical proofs.
\end{abstract}

\section{\label{sec: intro}Introduction}

The one-dimensional Kuramoto--Sivashinsky equation (KSE) has attracted significant attention not only as a physical model \cite{Kuramoto1975, Kuramoto1976, Sivashinsky1977, Michelson1977, Sivashinsky1980b} but also as a bridge between low-dimensional chaotic systems and high-dimensional spatiotemporal chaos, including fluid turbulence \cite{Hyman1986, Holmes1997, Holmes2012}. Here we consider mean-zero solutions $u(x,t)$ with spatial period $2\pi L$, where $L$ is a parameter:
\beq
\begin{aligned}
&u_t = 
-u\,u_x - u_{xx} - u_{xxxx}, \\[6pt]
&u(x,t)=u(x+2\pi L,t), \\[6pt]
&u(x,0) = u_0(x),\quad \int u_0(x)\,\dx=0.
\end{aligned}
\label{eq: ks}
\eeq
For simplicity we restrict most of our analysis to the subspace of odd solutions where $ u(-x,t)=-u(x,t)$. All spatial integrals are over $[-\pi L,\pi L]$ unless otherwise noted. With sufficiently smooth initial conditions, solutions exist for all $t\ge0$ and have strong regularity, including finiteness of all Sobolev norms \cite{Nicolaenko1985}. Simulations of the KSE display spatiotemporal chaos of increasingly large dimension as the spatial period is raised \cite{Pomeau1984, Hyman1986a}.

Our objective is twofold: to develop methods for bounding average quantities governed by nonlinear PDEs, and to produce novel results for the KSE using these methods. The methods follow an approach used commonly in the study of dynamics: constructing a functional that, by satisfying suitable inequalities, implies the desired result. In particular we construct \emph{auxiliary functionals}, which satisfy conditions implying bounds on time averages. These conditions are related to but not the same as the conditions on Lyapunov functionals, which can be used to bound instantaneous quantities or prove nonlinear stability. In this work the choice of auxiliary functional is optimized with computer assistance. This is done by replacing all relevant inequalities with constraints that require certain polynomials to be representable as sums of squares. The resulting convex optimization problem can be translated into a semidefinite program (SDP) and solved numerically. Such SDP-based methods have produced various mathematical statements about nonlinear ODEs, including verification of nonlinear stability \cite{Parrilo2000, Papachristodoulou2002, Anderson2015a}, estimation of attractor basins \cite{Henrion2014}, and bounds on time averages \cite{Chernyshenko2014, Fantuzzi2016, Goluskin2018}. Here we extend methods for bounding time averages to the PDE setting.

The quantity on which we focus is the spatially averaged energy in the KSE,
\beq
\E(t) := \fint u(x,t)^2\,\dx, 
\label{eq: energy}
\eeq
where $\fint$ denotes the spatial average $\frac{1}{2\pi L}\int_{-\pi L}^{\pi L}$. Past authors have discussed the energy in various terms, including $\E$, $\tfrac12\E$, and $\Vert u\Vert_2=\left(\int u^2\dx\right)^{1/2}$. Here we compute upper bounds on the spatiotemporal average of the energy,
\beq
\ov{\E}:=\limsup_{T\to\infty}\frac{1}{T}\int_0^T\E(t)\,{\rm d}t,
\label{eq: energy avg}
\eeq
for various values of $L$. The scaling of our bounds with increasing $L$ has implications for the long-standing challenge to analytically prove bounds that scale optimally when $L\gg1$.

Upper bounds on the energy have been pursued for several decades as part of a larger conjecture dating to the 1980s \cite{Pomeau1984, Wittenberg1999}, which asserts that various spatial averages scale intensively in large domains---that is, they are independent of $L$ when $L\gg1$. Intensive scaling would be confirmed for a quantity of interest if one could prove upper and lower bounds that are $O(1)$ for $L\gg1$. The trivial lower bound $\ov\E \geq 0$ cannot be improved upon for general initial conditions since it is saturated by the $u=0$ state. On the other hand, there is hope of constructing $O(1)$ upper bounds since no solutions displaying larger scaling have been reported. To date, $O(1)$ upper bounds have been proved only for steady solutions \cite{Michelson1986}, and logarithmic-in-$L$ upper bounds have been proved only for spatiotemporal averages of $(|\partial_x|^\alpha u)^2$ with $1/3\le\alpha\le2$ \cite{Otto2009, Goldman2015}. The best bounds proven for mean energy grow as powers of $L$.

Most bounds on energy in the KSE have focused on its instantaneous value at late times,
\beq
\E_\infty := \limsup_{t\to\infty}\E(t).
\label{eq: E_inf}
\eeq
Upper bounds on $\E_\infty$ apply to the time average $\ov\E$, but not vice versa, since $\ov\E\le\E_\infty$. The first result, proved in the 1980s, was that $\E_\infty\le O(L^4)$ for odd solutions \cite{Nicolaenko1985}. Later this was extended to all solutions of \eqref{eq: ks} and improved to $\E_\infty\le O(L^2)$ \cite{Collet1993, Goodman1994, Bronski2006}. All of these results were proved using quadratic Lyapunov functionals. Bounds that grow more slowly than $O(L^2)$ have been proved by Otto and coworkers using an ``entropy method'' in which solutions of the KSE are estimated using entropy solutions of the inviscid Burgers equation. Their bounds on the instantaneous \cite{Giacomelli2005} and time-averaged \cite{Otto2009, Goldman2015} energy are currently the best available in the $L\gg1$ limit:
\beq
\E_\infty\le o(L^2)
\qquad \text{and} \qquad
\ov{\E} \le O(L^{2/3+}).
\label{eq: best bounds}
\eeq
Here we focus on bounding $\ov\E$, rather than $\E_\infty$, for two reasons. First, it is easier in the time-averaged case to confirm that bounds are sharp because, it turns out, they are saturated by simple equilibria. Second, time-averaged bounds are easier to optimize computationally because the optimization problem is convex, as described in the next section.

Despite the success of the entropy method, which relies on particular features of the KSE, there are compelling reasons to continue developing the earlier approach. Lyapunov functionals and auxiliary functionals, which can be used to bound $\E_\infty$ and $\ov\E$, respectively, are very broadly applicable to nonlinear PDEs. Furthermore, our results for the KSE suggest a way to possibly improve upon \eqref{eq: best bounds}. It has been argued that $\E_\infty\le O(L^2)$ is the best bound provable using quadratic Lyapunov functionals \cite{Bronski2006, Wittenberg2014}. Indeed, using SDPs to optimize quadratic Lyapunov functionals improves the bounds of \cite{Bronski2006} but not their $O(L^2)$ scaling~\cite{Fantuzzi2015}. We likewise find for the time-averaged problem that quadratic auxiliary functionals give $\ov\E\le O(L^2)$. However, we also produce better bounds on $\ov\E$ using auxiliary functionals that are not quadratic. In particular, we use SDP-based computations to construct quartic and sextic polynomial auxiliary functionals. Since the restriction to quadratic functionals is equivalent to the ``background method'' used in past studies of the KSE and numerous fluid dynamical models \cite{Chernyshenko2017}, the present work is a generalization of the background method.

Our computational approaches require approximating a PDE by systems of ODEs that are derived from the PDE by Galerkin truncation. Methods for using SDP computations to bound time averages in ODE systems have been demonstrated previously \cite{Chernyshenko2014, Fantuzzi2016, Goluskin2018}, albeit only for systems of dimension 9 or fewer. One option is to approximate bounds for the PDE by computing bounds for ODE truncations, increasing the truncation dimension until the bounds converge. A second option is to compute bounds for the PDE using an ODE truncation of fixed size, along with analytical estimates controlling the difference between solutions to the ODE system and the PDE. The bound for the PDE is computed by solving an SDP that incorporates the truncated system and the analytical estimates. Here we apply both approaches to the KSE. Each approach produces bounds for various $L$ that are within 0.1\% of being sharp. The first approach does so with less computational cost, so most of our results have been computed in this way.

Section \ref{sec: method} describes how auxiliary functionals can be used to bound spatiotemporal averages in PDEs, and how the analysis can be carried out with computer assistance by solving SDPs. The relationship to Lyapunov functionals is explained, and the methodology is tailored to the KSE. Section \ref{sec: steady} reviews the bifurcation structure of the KSE as needed to judge the quality of our bounds, which are presented in  \S\ref{sec: results}. In  \S\ref{sec: con} we offer several conjectures about the KSE motivated by our computations, followed by conclusions about our methodology in general. Appendices provide two propositions that help constrain auxiliary functionals, details on computations, and a reinterpretation of some results in terms of the background method.

\section{\label{sec: method}Auxiliary functionals and Lyapunov functionals}

Bounds on time averages for many different PDEs have been proved using auxiliary functionals, although often the proofs are not presented in this framework, and the choice of functional is not explicit. The framework is summarized in \S\ref{sec: framework}, and the corresponding framework for bounding instantaneous values using Lyapunov functionals is discussed in \S\ref{sec: lyap}. As a simple example for the KSE, \S\ref{sec: quad} describes quadratic auxiliary and Lyapunov functionals, which are tantamount to the background method. Section \ref{sec: poly af} describes a general framework for polynomial auxiliary functionals of any degree. The construction of such functionals can be posed as a polynomial optimization problem that is equivalent to an SDP. This is described in \S\ref{sec: general trunc} for Galerkin truncations of PDEs and in \S\ref{sec: general pde} for the full PDEs. In \S\ref{sec: poly af kse} both approaches are formulated for the KSE in particular.

\subsection{\label{sec: framework}Auxiliary functionals}

Consider an autonomous PDE,
\beq
u_t(x,t)=F(u(x,t)),
\label{eq: dyn sys}
\eeq
whose righthand side $F(u)$ involves differentiation in $x\in\R^m$. Assume that solutions remain in some function space $\U$, meaning $u(\cdot,t)\in\U$ for all $t\ge0$. Suppose we are interested in a real-valued functional $\Phi(u(\cdot,t))$ that remains bounded along each solution trajectory. Its infinite-time average $\ov\Phi$ may depend on the trajectory along which it is evaluated. Our objective is to bound the possible values of $\ov\Phi$. This can be done using an \emph{auxiliary functional}, $V:\U\to\R$, that is absolutely continuous, differentiable, and remains bounded along trajectories. The conditions on $V$ ensure that the infinite-time average $\ov{\ddt V}$ vanishes, so along each solution $u(x,t)$ of the PDE,
\beq
\ov{\Phi} = \ov{\Phi + \tfrac{\rm d}{{\rm d}t}V}.
\label{eq: identity}
\eeq
Let $D_V:\U\to\R$ be a time-independent functional that is equal to $\ddt V(u(\cdot,t))$ along all trajectories. An expression for $D_V(u)$ must be deduced from the dynamics \eqref{eq: dyn sys}. For instance if $V(u)=\int u^2\dx$ and the PDE solution $u(x,t)$ is sufficiently regular, then $D_V(u)=2\int u\,F(u)\dx$. To prove an upper bound $\ov\Phi\le B$ it suffices to find a $V$ such that $\Phi(u) + D_V(u)\le B$ for all $u\in\U$. This condition implies $\Phi(u(\cdot,t)) + \ddt V(u(\cdot,t))\le B$ on all trajectories at all times, which implies the desired bound on $\ov\Phi$ via the identity \eqref{eq: identity}. For convenience we rearrange the sufficient condition for $\ov\Phi\le B$ into a standard form:
\beq
S(u) := B - \Phi(u) - D_V(u) \ge 0 \quad \forall u\in\U.
\label{eq: S}
\eeq
The PDE enters the functional $S(u)$ only through the derivation of $D_V(u)$.

If all auxiliary functionals $V$ in some class $\V$ are known to be admissible, meaning $\ov{\ddt V}$ vanishes on all PDE trajectories, we can optimize over $\V$ to minimize the bound:
\beq
\ov\Phi \le \inf_V B \quad \text{s.t.}\quad 
\begin{array}[t]{l}V\in\V, \\ 
S(u)\ge0\quad\forall u\in\U.\end{array}
\label{eq: V opt}
\eeq
If the set $\V$ is convex, the righthand side of \eqref{eq: V opt} is a convex optimization problem. This fact is central to our SDP-based methods for optimizing $V$ computationally. 

When properties of the optimal $V$ can be anticipated, one can restrict attention to such members of $\V$. Two properties that arise often, including in our application to the KSE, concern the symmetry and boundedness of $V$. First, if the governing equations and quantity to be bounded are invariant under some symmetry, then the optimal $V$ is also. This is made precise in the ODE case by Proposition~\ref{prop: sym} in Appendix~\ref{app: prop}. Second, if $|\Phi(u)|<\infty$ implies $|V(u)|<\infty$, and $\Phi$ is bounded on the global attractor, then only $V$ which are bounded below can give finite bounds on $\ov\Phi$. This is made precise by Proposition~\ref{prop: bounded} in Appendix~\ref{app: prop}.
 
For broad families of PDEs \cite{TobascoPC} and ODEs \cite{Tobasco2018}, there exists a convex class $\V$ such that the righthand side of \eqref{eq: V opt} is guaranteed to produce arbitrarily sharp bounds:
\beq
\sup_{u(x,t) \text{ solves \eqref{eq: dyn sys}}}\ov\Phi =
\inf_V B \quad \text{s.t.}\quad 
\begin{array}[t]{l}V\in\V, \\ 
S(u)\ge0\quad\forall u\in\U.
\end{array}
\label{eq: duality}
\eeq
In practice it may or may not be tractable to optimize over a $\V$ that is large enough to give a bound $B$ saturating, or nearly saturating, the above equality. If not, it still is desirable to optimize bounds over a numerically tractable subset of $\V$.

\subsection{\label{sec: lyap}Lyapunov functionals}

Often when studying a dynamical PDE \eqref{eq: dyn sys}, the objective is to bound instantaneous values $\Phi(u(\cdot,t))$ on a global attractor, rather than the time average $\ov\Phi$. Suppose we seek an upper bound on the long-time maximum defined on each trajectory $u(x,t)$ by
\beq
\Phi_\infty := \limsup_{t\to\infty}\Phi(u(\cdot,t)).
\label{eq: Phi_inf}
\eeq
An upper bound $\Phi_\infty\le C$ for all trajectories can be proved by finding some functional $\LL:\U\to\R$, and the corresponding functional $D_W:\U\to\R$ that is equal to $\ddt W(u(\cdot,t))$ along all PDE trajectories, such that
\begin{align}
\Phi(u) &\le \LL(u), \label{eq: lyap 1} \\
a\,D_W(u) &\le C - \LL(u) \label{eq: lyap 2}
\end{align}
for some $a>0$. We regard $W(u)$ as a kind of Lyapunov functional, although in some cases it can be negative. The above conditions are not the only ones that can be used to bound $\Phi_\infty$, but they are among the simplest. Condition~\eqref{eq: lyap 2} implies $\LL_\infty\le C$ by Gronwall's inequality, provided that $\ddt W(u(\cdot,t))$ is continuous along all trajectories, and then~\eqref{eq: lyap 1} implies $\Phi_\infty\le C$.

For any Lyapunov functional $\LL$ that proves $\Phi_\infty\le C$ via conditions \eqref{eq: lyap 1}--\eqref{eq: lyap 2}, the auxiliary functional $V(u)=a\,W(u)$ proves the same upper bound $\ov\Phi\le C$ via condition \eqref{eq: S}. The converse is not true, which is consistent with the fact that $\ov\Phi\le\Phi_\infty$.

As with the optimization of time-averaged bounds in \eqref{eq: V opt}, it is natural to seek the smallest bound $C$ by optimizing $\LL$ within some convex class, subject to \eqref{eq: lyap 1}--\eqref{eq: lyap 2}. This optimization is convex for fixed $a$, but the joint optimization over $a$ and $W$ is not convex because $a$ and $D_W$ are multiplied in condition \eqref{eq: lyap 2}. Optimizing $W$ computationally while sweeping through $a$ values works well for chaotic ODEs \cite{Goluskin2018b}, and an extension to PDEs is ongoing, but here we focus on the easier optimization \eqref{eq: V opt} over bounds on time averages.

\subsection{\label{sec: quad}Quadratic auxiliary functionals for the KSE}

For odd solutions of the KSE, the simplest auxiliary functionals giving finite bounds on $\ov\E$ are quadratic functionals in the class
\beq
\V=\left\{ c \int u^2\dx+ P(u) ~\big|~c\in\R,~ P:\U\to \R\right\},
\label{eq: V quad}
\eeq
where $P$ is a linear functional. All $V\in\V$ remain bounded along mean-zero solutions of the KSE, so it remains only to enforce $S(u)\ge0$. This requires $D_V(u)$ to be bounded above, which is why the quadratic term in $V$ is constrained to be as in \eqref{eq: V quad}. A quadratic term in $V$ generically leads to a cubic term in $D_V$ since the nonlinearity of the KSE is quadratic. This is not the case for the quadratic term of \eqref{eq: V quad} since the nonlinearity of the KSE conserves energy:
\beq
\frac{{\rm d}}{{\rm d} t}\frac12\int u^2\,\dx = \int u u_t\,\dx = -\int u(uu_x+u_{xx}+u_{xxxx}) \dx=\int(u_x^2-u_{xx}^2) \dx.
\eeq

One way to ensure that $V$ is in the class \eqref{eq: V quad} is to use the ansatz 
\beq
V(u) = \frac{\alpha}{2}\fint (u-\zeta)^2 \dx.
\label{eq: V background}
\eeq
The quadratic terms in \eqref{eq: V quad} and \eqref{eq: V background} are identical for $\alpha=4\pi L c$, and the $u$-independent term in \eqref{eq: V background} is irrelevant since it does not enter $\ddt V$. As for the linear term, any $\zeta$ for which $\int u\zeta \dx$ is well defined gives a linear functional $P$. (Whether any such $P$ can be represented by a corresponding $\zeta$ depends on the class $\U$.) Using the ansatz \eqref{eq: V background} is an instance of the background method \cite{Doering1992}, so called because $u-\zeta$ is the deviation of $u$ from a ``background'' function $\zeta(x)$.

With the $V$ ansatz \eqref{eq: V background} used to determine $D_V$, the functional $S$ defined by \eqref{eq: S} is
\beq
S(u) = 
\fint\left[ \alpha u_{xx}^2-\alpha u_x^2+\left(\tfrac12\alpha\zeta_x-1\right)u^2 \right] \dx
+ 
\fint \left(\alpha\zeta_xu_x-\alpha\zeta_{xx}u_{xx}\right) \dx + B.
\label{eq: S quad}
\eeq
We regard the tunable variables as $\alpha$ and $\alpha\hspace{1pt}\zeta$, as opposed to $\alpha$ and $\zeta$, since $S$ is jointly convex in the former pair. For quadratic $V$, the optimization problem \eqref{eq: V opt} for the best upper bound on $\ov\E$ in the KSE is
\beq\
\ov\E \le \inf_{\alpha,\,\alpha\zeta} B \quad \text{s.t.}\quad
\begin{array}[t]{l}
\alpha\in\R \\ \zeta \in C^2 \\ 
\eqref{eq: S quad}\ge0 \quad\forall u\in C^2.
\end{array}
\label{eq: V opt quad}
\eeq
Numerical solutions of the righthand minimization are $O(L^2)$ for $L\gg1$, as reported in \S\ref{sec: results}, and Appendix~\ref{app: background} gives examples of optimal $\zeta(x)$ for several domain sizes. If upper bounds smaller than $O(L^2)$ are to be proved using auxiliary functionals, we must generalize beyond the quadratic class \eqref{eq: V quad}. 

\begin{remark}
Like auxiliary functionals, every quadratic Lyapunov functional must have a leading term proportional to $\int u^2\dx$ in order for condition \eqref{eq: lyap 2} to be possible. Thus any quadratic Lyapunov functional can be expressed as in \eqref{eq: V background}, up to addition of another constant. This is the framework that past authors have used when proving bounds on $\E_\infty$ using Lyapunov functionals \cite{Nicolaenko1985, Collet1993, Goodman1994, Bronski2006, Fantuzzi2015}. Successively better analytical constructions for $\zeta$ led to the $O(L^2)$ bounds of \cite{Bronski2006}. Still better $\zeta$ computed using SDPs improved these bounds but not their scaling \cite{Fantuzzi2015}. The full optimization over quadratic Lyapunov functionals subject to \eqref{eq: lyap 1}--\eqref{eq: lyap 2} has not been solved, but we suspect that it would produce $O(L^2)$ bounds also. Optimizing over higher-degree functionals to bound $\E_\infty$, as we do here for the time-averaged problem, remains a topic for future work.
\end{remark}

\begin{remark}
Most applications of the background method---that is, of quadratic auxiliary functionals with leading terms proportional to energy---have been to the Navier--Stokes equations and related fluid dynamical systems (e.g., \cite{Doering1994, Doering1996, Nicodemus1997, Plasting2003, Goluskin2016a}). These PDEs have the same type of quadratic nonlinearity as the KSE. Thus, as with the KSE, any term of $V$ that is quadratic in the velocity field $\mathbf{u}(\mathbf{x},t)$ must be proportional to the energy in order for $D_V$ to be quadratic instead of cubic. Therefore such $V$ can be written as
\beq
V(\mathbf{u})=c\int_\Omega|\mathbf{u}-\boldsymbol{\zeta}|^2{\rm d}\mathbf{x},
\eeq
possibly with additional terms involving other variables such as temperature. As pointed out by Chernyshenko \cite{Chernyshenko2017}, such quadratic auxiliary functionals underlie all past uses of the background method to bound time averages, although in many cases $V$ was not given explicitly. 
\end{remark}

\subsection{\label{sec: poly af}Auxiliary functionals of any polynomial degree}

Having illustrated the use of quadratic auxiliary functionals for the KSE, let us return to the general PDE \eqref{eq: dyn sys} and the optimization \eqref{eq: V opt} over auxiliary functionals. Our aim is to compute auxiliary functionals of higher polynomial degree by methods of polynomial optimization. In order for polynomial optimization to be directly applicable, we restrict attention to dynamics and quantities to bound---$F(u)$ and $\Phi(u)$---with polynomial dependence on $u$. That is, scaling $u\mapsto\beta u$ would rescale each term in $F$ and $\Phi$ by a nonnegative integer power of $\beta$. The KSE and its mean energy have this property. We also restrict attention to auxiliary functionals with polynomial dependence on $u$. Whereas quadratic functionals of the restricted form \eqref{eq: V quad} can be represented by a constant $c$ and background function $\zeta$, representations of higher-degree $V$ are more complicated. Here we let $V$ be a polynomial function of the projections of $u$ onto an $L^2$-orthonormal basis $\{u_n(x)\}_{n\geq 1}$. That is, we expand
\beq
u(x,t) = \sum_{n=1}^\infty a_n(t)u_n(x)
\label{eq: Galerkin inf}
\eeq
and consider
\beq
V(u) = V(a_1,a_2,\ldots),
\label{eq: V inf poly}
\eeq
where the righthand side is a polynomial function of the expansion coefficients $\{a_n\}_{n\geq 1}$.

For analytical purposes $V$ might include an infinite number of monomials, but for computational purposes we must optimize $V$ over a finite-dimensional space. A relatively simple approach is to truncate the Galerkin expansion~\eqref{eq: Galerkin inf} after $N$ modes, replace the PDE with its projection onto these modes, and let $V=V(a_1,\ldots,a_N)$. This is described in \S\ref{sec: general trunc}. An alternate approach, which does not require truncating the PDE, is to let $V$ depend not only on $(a_1,\,\ldots,\,a_N)$ but also on functionals of the ``tail'' of the expansion of $u$. This is described in \S\ref{sec: general pde}.

\subsubsection{\label{sec: general trunc}Bounds for truncated ODE systems}

To approximate the PDE \eqref{eq: dyn sys} by a finite-dimensional system we employ an $N$-mode Galerkin projection. Truncating the expansion \eqref{eq: Galerkin inf} after the first $N$ terms and integrating each basis function $u_n(x)$ against the PDE yields a system of ODEs of the form
\beq
\ddt\ba = \f(\ba),
\label{eq: ode}
\eeq
where  $\ba=(a_1,\ldots,a_ N)$ is the vector of mode amplitudes, and each component of $\f(\ba)$ is a polynomial since the PDE has polynomial dependence on $u$. We likewise replace the quantity to be bounded with its $N$-mode truncation, $\Phi_N(\ba)$, and seek bounds on the time average~$\ov\Phi_N$.

If all trajectories of the truncated system \eqref{eq: ode} remain bounded, then bounds on $\ov\Phi_N$ can be computed using existing SDP-based methods for ODEs \cite{Chernyshenko2014, Fantuzzi2016, Goluskin2018}. The chain rule gives $\ddt V(\ba)=\f\cdot\nabla V(\ba)$, so the bounding condition \eqref{eq: S}, which proves $\ov\Phi_N\le B$, is the polynomial inequality
\beq
S(\ba) = B - \Phi_N(\ba) - \f\cdot\nabla V(\ba) \ge 0 \quad \forall\ba\in\R^N.
\label{eq: S N}
\eeq
Even with computer assistance it is prohibitively difficult in general to decide whether a polynomial is nonnegative---the computational complexity is NP-hard \cite{Murty1987}. We thus employ a standard technique in polynomial optimization \cite{Parrilo2013a}: replace the condition $S(\ba)\ge0$ with the stronger condition that $S$ can be represented as a sum of squares of other polynomials. 
Precisely, let $\R[\ba]_{N,2d}$ denote the set of degree-$2d$ polynomials in $\ba\in\R^N$, and let $\Sigma[\ba]_{N,2d}$ denote its subset of sum-of-squares (SOS) polynomials,
\begin{equation}
\Sigma[\ba]_{N,2d} := \left\{
p\in\R[\ba]_{N,2d} \,:\, 
\exists k\in\mathbb{N},\,
q_1,\,\ldots,\,q_k \in\R[\ba]_{N,d} \text{ such that } 
p(\ba) = \sum_{i=1}^k q_i(\ba)^2
\right\}.
\end{equation}
We enforce $S(\ba)\ge0$ by requiring $S\in\Sigma[\ba]_{N,2d}$. The degree of $V$ can be no larger than that of $S$ and may need to be smaller, so it suffices to choose $V$ from $\R[\ba]_{N,2d}$. With the degree of $S$ no greater than $2d$, the best bound that can be proved using the SOS framework is
\beq
\ov\Phi_N \le B_{N,2d} := \min_V B \quad \text{s.t.}\quad 
\begin{array}[t]{l}V\in\R[\ba]_{N,2d}, \\ 
S\in\Sigma[\ba]_{N,2d}.\end{array}
\label{eq: B opt trunc}
\eeq
The tunable parameters are the coefficients in $V$, which appear linearly in $V$ and $S$, and the bound $B$, which appears linearly in $S$. Thus the SOS constraint is convex, and the minimization problem can be formulated as an SDP and solved computationally~\cite{Parrilo2013a}.

The bounds on truncated averages $\ov\Phi_N$ found by solving \eqref{eq: B opt trunc} are not necessarily bounds on $\ov\Phi$ in the full PDE, so we consider the limit of these bounds:
\beq
B_{2d} := \limsup_{N\to\infty} B_{N,2d}.
\label{eq: B_2d}
\eeq
The bound $\ov\Phi\le B_{2d}$ holds for the PDE if for every PDE solution $u(x,t)$ there exists a sequence of solutions to truncated systems where $\ov\Phi_N$ converges to $\ov\Phi$. This appears to be true for $\ov\E$ in the KSE, but we do not prove it here. In such cases, the PDE bound $B_{2d}$ can be approximated by increasing $N$ until $B_{N,2d}$ converges numerically, much like choosing spatial resolution when numerically integrating PDEs. In the computations for the KSE reported here the numerical values of $B_{N,2d}$ converge quickly. Nonetheless, since one obtains only numerical approximations of $B_{2d}$, these values are not guaranteed to be valid bounds for the PDE. For many purposes this is not important, and the approach described above would be practical. For other purposes, such as computer-assisted proofs, it would be preferable to obtain bounds for the full PDE at finite $N$. This motivates the alternative framework of the following subsection, where the deviation between the truncated system \eqref{eq: ode} and the full PDE \eqref{eq: dyn sys} is estimated rigorously.

\subsubsection{\label{sec: general pde}Bounds for the full PDE}

To construct an auxiliary functional for the full PDE \eqref{eq: dyn sys} we cannot simply truncate the Galerkin expansion \eqref{eq: Galerkin inf} after $N$ terms. Instead we retain the tail $v(x,t)$ of the expansion,
\beq
u(x,t) = \sum_{n=1}^Na_n(t)u_n(x) + v(x,t).
\label{eq: Galerkin}
\eeq
Largely following the ideas of Goulart and Chernyshenko \cite{Goulart2012}, we let the polynomial $V$ depend not only on $\ba$ but also on $\bb=(b_1,\ldots,b_M)$, where each functional $b_m(v(\cdot,t))$ depends only on the tail. Assume the functionals included in $\bb$ allow $\Phi(u)$ to be represented exactly as $\Phi(\ba,\bb)$. Expressions for the evolution of $\ba$ and $\bb$ can be derived from the PDE in the general form
\begin{align}
\ddt\ba &= \f(\ba) + \mathbf G(\ba,v), \label{eq: uncert 1 gen} \\
\ddt\bb &= \mathbf H(\ba,v), \label{eq: uncert 2 gen}
\end{align}
where $\f(\ba)$ is the ODE truncation \eqref{eq: ode}. The functional $S$ whose nonnegativity would prove $\ov\Phi\le B$ is then
\beq
S(u) = B - \Phi(\ba,\bb) - \f(\ba)\cdot\nabla_\ba V(\ba,\bb) - 
	\left[ \mathbf G(\ba,v)\cdot\nabla_\ba V(\ba,\bb) 
		+ \mathbf H(\ba,v)\cdot\nabla_\bb V(\ba,\bb) \right]
\eeq
The expression inside the square brackets need not be determined by $(\ba,\bb)$ alone since it can depend more generally on the tail $v$. Likewise equations \eqref{eq: uncert 1 gen}--\eqref{eq: uncert 2 gen} may not fully determine the evolution of $(\ba,\bb)$, in which case they form an ``uncertain system'' \cite{Goulart2012}. However, by estimating the bracketed expression analytically in terms of $(\ba,\bb)$, possibly with the help of SOS polynomial constraints, one can seek a polynomial $T(\ba,\bb)$ such that $S(u)\ge T(\ba,\bb)$ for all $u\in\U$. Then the bounding condition $S(u)\ge0$ can be replaced by the stronger condition that $T(\ba,\bb)$ is an SOS polynomial, and the optimization over $V(\ba,\bb)$ can be formulated as an SDP. This approach is illustrated for the KSE in \S\ref{sec: non-quad pde}, where the only functional needed in $\bb$ is the energy of $v$.

\subsection{\label{sec: poly af kse}Higher-degree auxiliary functionals for the KSE}

Polynomial auxiliary functionals for the KSE have many degrees of freedom in their lower-order terms, but their leading terms are highly constrained. Proposition~\ref{prop: bounded} in Appendix~\ref{app: prop} implies that because $\E(u)$ is bounded below and remains finite along all trajectories, only $V(u)$ that are bounded below can give finite bounds on $\ov\E$. Since $V$ and $D_V$ must be bounded below and above, respectively, the leading terms of both functionals must be of even degree. Thus we constrain the leading term of $V$ as in
\beq
V(u) = c\left(\int u^2\dx\right)^d + P(a_1,a_2,\ldots),
\label{eq: V general}
\eeq
where $d\ge1$ is an integer and $P$ is a polynomial of degree at most $2d-1$. The leading term of $V$ is conserved by the nonlinearity of the KSE, so $D_V$ has a leading term of degree $2d$, rather than $2d+1$. There exist other degree-$2d$ expressions that are conserved by the nonlinearity, such as $\int u^{2d}\dx$, but we have obtained finite bounds only with the leading term in \eqref{eq: V general}. Similarly, the non-quadratic Lyapunov functionals constructed for the Navier--Stokes equations in \cite{Huang2015} have powers of the energy as their leading terms.

Since $u$ is odd and periodic we use the orthonormal sine basis,
\beq
u_n(x)=(\pi L)^{-1/2}\sin(n x/L),
\label{eq: basis}
\eeq
for the Galerkin expansion \eqref{eq: Galerkin inf}. The projections $a_n(t)$ are Fourier coefficients of $u(x,t)$, up to rescaling by $\sqrt{\pi L}$. We let $P$ in~\eqref{eq: V general} depend on the first $N$ projections, and possibly on the energy of the tail also. All $V$ defined in this way remain bounded along solutions of the KSE, as required.

Symmetry provides an additional constraint on $V$. The subspace of odd periodic solutions of the KSE is invariant under $u(x,t)\mapsto -u(x+\pi L,t)$, as is the energy $\E$. Thus it suffices to consider $V$ with this same symmetry, as implied by Proposition~\ref{prop: sym} in Appendix~\ref{app: prop}. In terms of Fourier coefficients, such $V$ are invariant under $a_n\mapsto (-1)^na_n$. Therefore, the polynomial $P$ in \eqref{eq: V general} can omit all monomials $a_{i_1}a_{i_2}\cdots a_{i_K}$ (with possible repetition in subscript values) where the sum of subscript values is odd. For instance, $V$ of quartic degree might contain the monomial $a_3^2a_6$ but not $a_1a_2a_4$. In the case of quadratic $V$, this symmetry corresponds to a background function $\zeta(x)$ that is $\pi L$-periodic, as illustrated in Appendix~\ref{app: background}.

\subsubsection{\label{sec: non-quad trunc}Truncated KSE}

In the sine basis, the $N$-mode Galerkin truncation \eqref{eq: ode} of the KSE is \cite{Papageorgiou1991}
\begin{equation}
f_n(\ba) = \left(\frac{n}{L}\right)^2\left[1-\left(\frac{n}{L}\right)^2\right]a_n 
	+ \frac1{\sqrt{\pi L}}\frac{n}{L}\left[ \frac12 \sum_{m=1}^{N-n}a_ma_{m+n}
	-\frac14\sum_{m=1}^{n-1}a_ma_{n-m}\right].
\label{eq: f}
\end{equation}
If bounds are to be computed only for truncated systems, many other bases could be chosen, although it would be more work to derive $f_n(\ba)$. To compute bounds for the full PDE as described in the next subsection, the sine basis has additional advantages due to each $u_n(x)$ being an eigenfunction of the KSE's linear operator.

The truncation of spatially averaged energy is $\E_N=\tfrac{1}{2\pi L} |\ba|^2$, and the truncation of the auxiliary functional \eqref{eq: V general} is
\beq
V(\ba) = c|\ba|^{2d} + P(\ba).
\label{eq: V trunc}
\eeq
The polynomial $P(\ba)$ in $N$ variables has degree $2d-1$ or less. That is, $P\in\R[\ba]_{N,2d-1}$. The best bound that can be proved using the optimization \eqref{eq: B opt trunc} with $V(\ba)$ of degree $2d$ is
\beq
\ov\E_N \le B_{N,2d} := \min_{c,\,P} B \quad \text{s.t.} \quad 
\begin{array}[t]{l} c\in\R, \\ 
P\in \R[\ba]_{N,2d-1}, \\
S\in\Sigma[\ba]_{N,2d},
\end{array}
\label{eq: B trunc kse}
\eeq
where $S$ is defined as in \eqref{eq: S N} with $\f$ given by \eqref{eq: f}. We approximate bounds for the full KSE using these $B_{N,2d}$ by increasing $N$ to approach the large-$N$ limit $B_{2d}$ defined in \eqref{eq: B_2d}. In the case of quadratic $V$, the limit $B_2$ is the optimal bound \eqref{eq: V opt quad} provable by the background method, and the coefficients of the linear function $P(\ba)$ are proportional to the Fourier coefficients of the background function $\zeta(x)$ as described in Appendix~\ref{app: background}. Increasing the polynomial degree $2d$ can only improve bounds or leave them unchanged. Here we have computed numerical approximations to $B_{2d}$ over a range of domain sizes $L$ for auxiliary functionals of degree $2d=2$, 4, and 6.

\subsubsection{\label{sec: non-quad pde}Full KSE}

To bound energy in the KSE using the framework of \S\ref{sec: general pde}, it suffices to let the auxiliary functional \eqref{eq: V general} depend on the first $N$ sine mode amplitudes, $\ba=(a_1,\ldots,a_N)$, and the energy of the tail,\footnote{Here $q^2$ is twice as large as the quantity of the same name in \cite{Goulart2012}.} 
\beq
q^2(t):=\int v(x,t)^2\dx = \sum_{n=N+1}^\infty a_n(t)^2.
\eeq
Such $V$ take the form
\beq
V(\ba,q^2)  =  c\left(|\ba|^2+q^2\right)^d+P(\ba,q^2),
\label{eq: V pde}
\eeq
where $P$ is a polynomial in the $N+1$ variables $(\ba,q)$ whose degree is no larger than $2d-1$. That is, $P\in\R[\ba,q]_{N+1,2d-1}$. The uncertain system \eqref{eq: uncert 1 gen}--\eqref{eq: uncert 2 gen} takes a particular form that arises also in the study of the Navier--Stokes equations \cite{Goulart2012}:
\begin{align}
\ddt\ba &= \f(\ba) + \bTh(\ba,v), \label{eq: uncert 1 kse} \\
\ddt\left(\tfrac12q^2\right) &= -\ba\cdot\bTh(\ba,v) + \Gamma(v), \label{eq: uncert 2 kse}
\end{align}
where here $\f$ is given by \eqref{eq: f} and
\begin{align}
\Gamma(v) &= \int v(-v_{xx}-v_{xxxx})\dx, \\
\Theta_n(\ba,v) &=  \frac{n}{2\pi^{1/2}L^{3/2}} \sum_{m=N-n+1}^\infty a_ma_{m+n}.
\end{align}
The functional \eqref{eq: S} that must be nonnegative is therefore
\beq
S(u) = B - \E -\f\cdot\nabla_\ba V - \bTh\cdot\bM - 2\dVdq\Gamma,
\label{eq: S pde}
\eeq
where 
\beq
\bM(\ba,q^2) = \nabla_\ba V -2\dVdq\ba = \nabla_\ba P -2\dPdq\ba.
\eeq
Although $\Gamma$ and $\bTh$ depend on the tail $v$, we can bound them analytically in terms of only $(\ba,q^2)$. This lets us derive a lower bound $S(u)\ge T(\ba,q^2)$, where $T$ is a polynomial.

It is simple to estimate $\Gamma$ in terms of $q^2$ since the sine basis is the eigenbasis of the linear operator $-(\partial_x^2+\partial_x^4)$. The corresponding eigenvalues are $\lambda_n=(n/L)^2-(n/L)^4$, so
\beq
\Gamma(v) \le \lambda_{N+1}q^2.
\label{eq: Gamma bound}
\eeq
Here we always choose $N>L$, in which case $\lambda_{N+1}<0$. In order for the above estimate to produce a lower bound on the last term in \eqref{eq: S pde}, we require that $\dVdq\ge0$ also. It is a convenient feature of the KSE that the Galerkin basis enabling the estimate \eqref{eq: Gamma bound} is the sine basis. In the case of the Navier--Stokes equations, obtaining an estimate analogous to \eqref{eq: Gamma bound} requires using the eigenbasis of the energy stability operator \cite{Goulart2012}, which often must be computed numerically.

To estimate $\bTh\cdot\bM$ in \eqref{eq: S pde} we apply the triangle inequality, $|\bTh\cdot\bM|\le \sum_{n=1}^N|\Theta_n||M_n|$, and then bound each $|\Theta_n|$ using the Cauchy--Schwarz inequality and Young's inequality,
\begin{align}
|\Theta_n| &\le \frac{n}{2\pi^{1/2}L^{3/2}}
	\left( \sum_{m=N+1}^\infty a_m^2 \right)^{1/2}
	\left( \sum_{m=N-n+1}^\infty a_m^2 \right)^{1/2} \notag \\
&\le \frac{n}{4\pi^{1/2}L^{3/2}}
	\left( 2q^2 + \sum_{m=N-n+1}^N a_m^2 \right). \label{eq: Theta bound}
\end{align}
To eventually obtain a polynomial lower bound on $S(u)$, we must bound each $|M_n|$ by an expression without absolute values. For this we introduce polynomials $R_n(\ba,q^2)$ and ensure that $|M_n|\le R_n$ for all $(\ba,q^2)$ by requiring
\beq
-R_n(\ba,q^2) \le M_n(\ba,q^2) \le R_n(\ba,q^2).
\label{eq: Rn}
\eeq
We find it suffices for $R_n$ to have the same polynomial degree as $M_n$. The preceding approach to bounding $\bTh\cdot\bM$ differs from that of \cite{Goulart2012}, where the term analogous to $\bTh\cdot\bM$ is bounded by $|\bTh\cdot\bM|\le\Vert\bTh\Vert_2\Vert\bM\Vert_2$, and then $\Vert\bTh\Vert_2$ is estimated. We have implemented both approaches for the KSE and find that using the sharper estimate \eqref{eq: Theta bound} greatly improves the eventual bounds on $\ov\E$.

A lower bound on the $S(u)$ functional \eqref{eq: S pde} follows from the above estimates on $\Gamma$, $|\Theta_n|$, and $|M_n|$. In particular, $S(u)\ge T(\ba,q^2)$ with
\begin{multline}
T(\ba,q^2) = B - \frac{1}{2\pi L}\left(|\ba|^2+q^2\right) - \f\cdot\nabla_\ba V - 2\lambda_{N+1}q^2\frac{\partial V}{\partial q^2} \\
- \frac{1}{4\pi^{1/2}L^{3/2}} \sum_{n=1}^NnR_n\left(2q^2 + \sum_{m=N-n+1}^N a_m^2 \right),
\end{multline}
provided that $\dVdq\ge0$ and each $R_n$ satisfies \eqref{eq: Rn}. Replacing each polynomial inequality with an SOS constraint gives:
\begin{align}
T(\ba,q^2)&\in\Sigma[\ba,q]_{N+1,2d} \label{eq: T SOS}, \\[2pt]
\dVdq(\ba,q^2)&\in\Sigma[\ba,q]_{N+1,2d-2}, \\[2pt]
R_n(\ba,q^2)-M_n(\ba,q^2) &\in\Sigma[\ba,q]_{N+1,2d-2}, 
	\quad \forall n\in\{1,\ldots,N\}, \\[2pt]
R_n(\ba,q^2)+M_n(\ba,q^2) &\in\Sigma[\ba,q]_{N+1,2d-2}, 
	\quad \forall n\in\{1,\ldots,N\}.
	\label{eq: R SOS 2}
\end{align}
The best bound that can be proved in this framework for $S$ of degree no more than $2d$ is
\beq
\ov\E\le B^{pde}_{N,2d} := \min_{c,\,P} B \quad \text{s.t.}\quad 
\begin{array}[t]{l}c\in\R,\\ 
P\in\R[\ba,q]_{N+1,2d-1},\\
\text{\eqref{eq: T SOS}--\eqref{eq: R SOS 2}.}
\end{array}
\label{eq: B pde kse}
\eeq
If $V(\ba,q^2)$ proves $\ov\E\le B^{pde}_{N,2d}$ for the full KSE in the above framework, $V(\ba,0)$ proves $\ov\E_N\le B_{N,2d}$ in the truncated framework \eqref{eq: B trunc kse} for some $B_{N,2d}\le B^{pde}_{N,2d}$. However, only the larger value $B^{pde}_{N,2d}$ is guaranteed to be a bound for the full KSE at finite $N$. We have solved \eqref{eq: B pde kse} computationally by converting it into an SDP. Results are reported in \S\ref{sec: results} for degree-4 auxiliary functionals over a range of domain sizes $L$. The resulting bounds converge from above as $N$ is raised.

\section{\label{sec: steady}Steady states of the KSE and its truncations}

\begin{figure}[t]
\begin{center}
\begin{tikzpicture}
\node at (0,0) {\includegraphics[height=140pt,trim={0 0 0 0},clip]{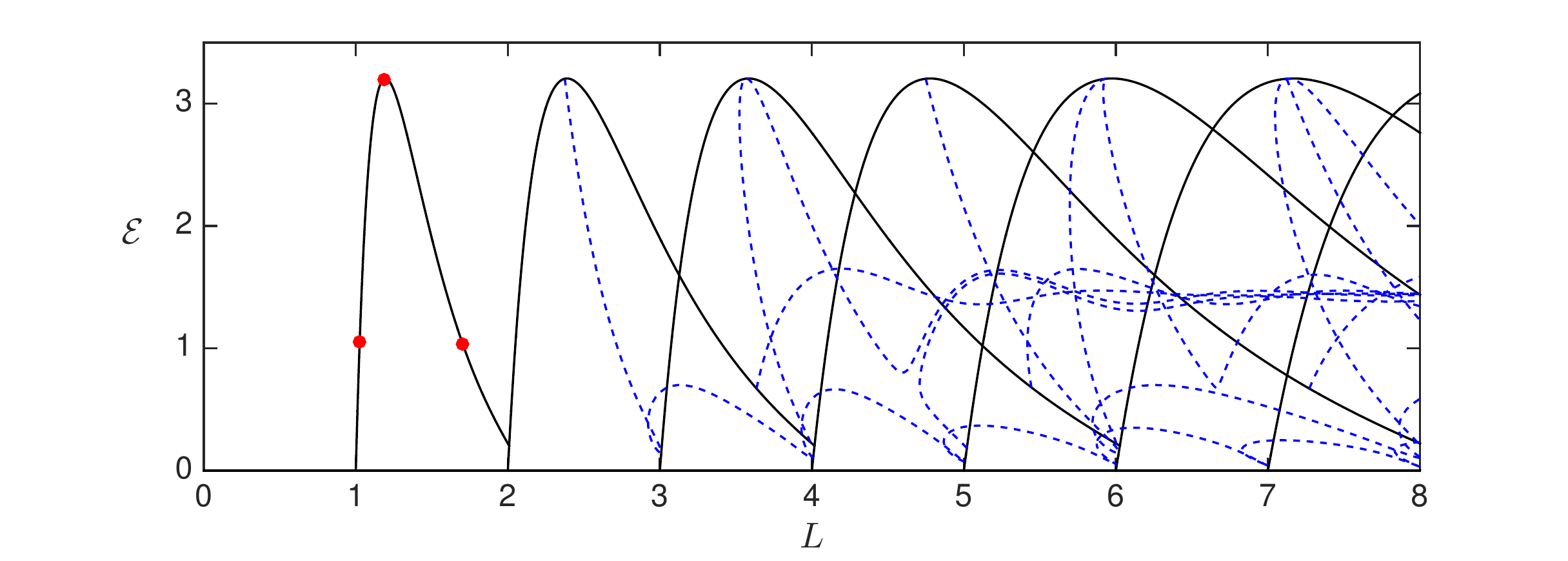}};
\node at (-4.55,1.6) {(a)};
\node at (-3.15,1.75) {{\footnotesize $E_1$}};
\node at (-1.45,1.75) {{\footnotesize $E_2$}};
\node at (.25,1.75) {{\footnotesize $E_3$}};
\node at (1.93,1.75) {{\footnotesize $E_4$}};
\node at (3.62,1.75) {{\footnotesize $E_5$}};
\node at (-3.96,-.55) {\color{red}\small\it b};
\node at (-3.75,1.76) {\color{red}\small\it c};
\node at (-3.07,-.55) {\color{red}\small\it d};
\end{tikzpicture}
\\
\begin{tikzpicture}
\node at (0,0) {\includegraphics[height=110pt,trim={15pt 5pt 2pt 8pt},clip]{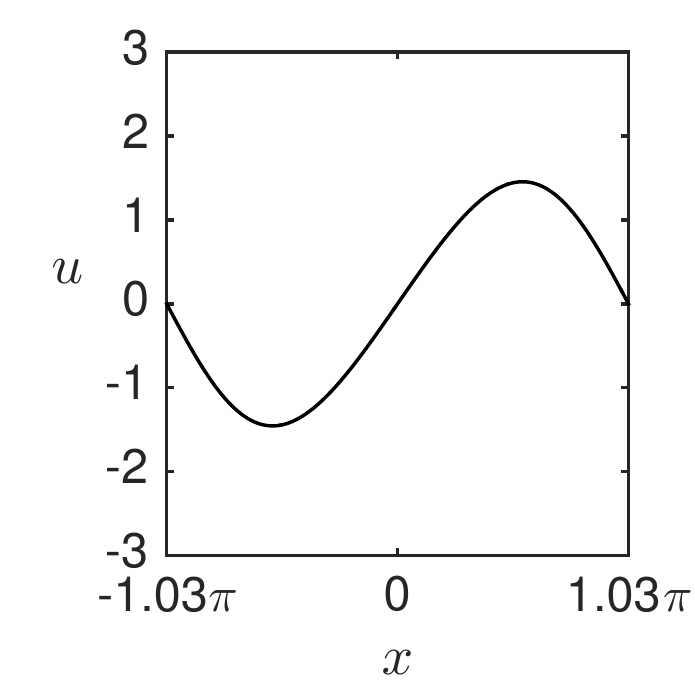}};
\node at (-.75,1.4) {(b)};
\end{tikzpicture}
\begin{tikzpicture}
\node at (0,0) {\includegraphics[height=110pt,trim={25pt 5pt 5pt 8pt},clip]{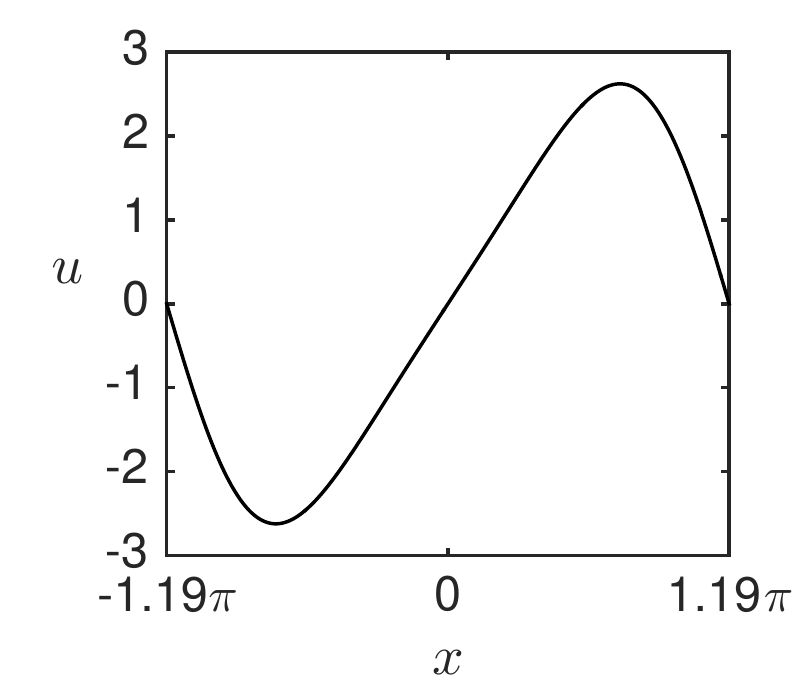}};
\node at (-1.15,1.4) {(c)};
\end{tikzpicture}
\begin{tikzpicture}
\node at (0,0) {\includegraphics[height=110pt,trim={25pt 5pt 15pt 8pt},clip]{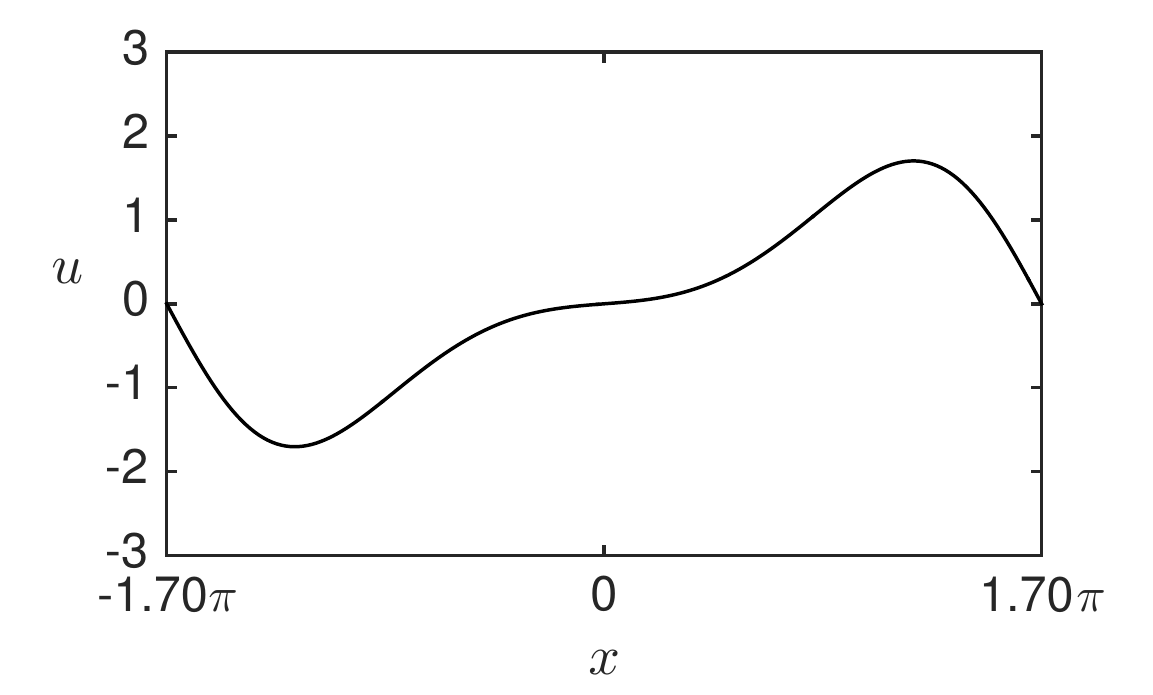}};
\node at (-2.1,1.4) {(d)};
\end{tikzpicture}
\end{center}
\caption{\label{fig: bif}(a) Bifurcation diagram showing the spatially averaged energy $\E$ of primary equilibria $E_n$ (\solidrule) and odd secondary equilibria ({\color{blue}\dashedrule}) of the KSE. Stability of equilibria is not indicated. (b-d) Steady solutions $u(x)$ on the $E_1$ branch at the three points indicated ({\footnotesize\color{red}$\bullet$}) in panel (a), where the values of $L$ are (b) 1.03, (c) 1.19, and (d) 1.70.}
% On 50<L<100, fit is  U~1.203*L^2.027.
\end{figure}
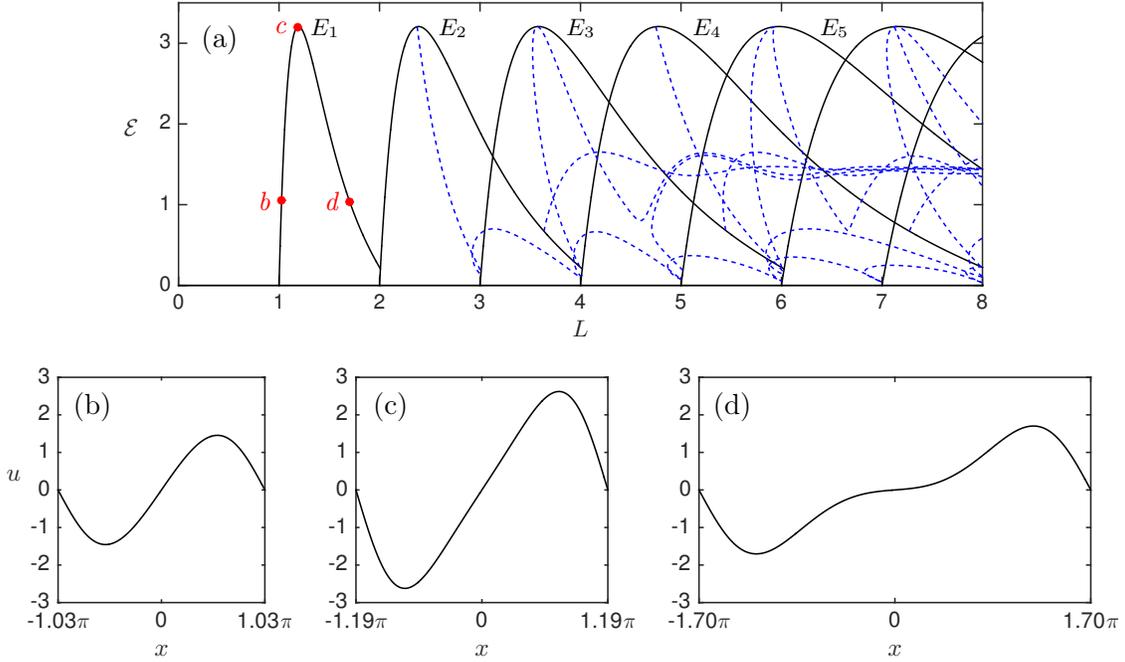

In order to judge the quality of the bounds reported in the next section, let us review the simplest odd steady states of the KSE \cite{Greene1988}. The zero state $u=0$ is globally attracting when $L<1$. As the domain size increases through $L=1$, the zero state becomes linearly unstable, and a bifurcation gives rise to four symmetry-related steady states that we call $E_1$. These states are mapped to one another by negation and/or translation by $\pi L$, so they have the same energy. The mean energy along the $E_1$ branch is shown in the bifurcation diagram of figure~\ref{fig: bif}(a), and $u(x)$ on the $E_1$ branch is plotted in figure~\ref{fig: bif}(b-d) for three different domain sizes. Rescaled versions of the bifurcation at $L=1$ occur at all integer values of $L$. This can be anticipated because for any solution of the KSE in a domain of size $2\pi L$, taking $n$ periods of this solution gives a solution in a domain of size $2\pi L n$. Thus for all integers $n>1$ there exist branches of ``primary equilibria'' $E_n$ consisting of $n$ copies of the $(2\pi L/n)$-periodic basic solution $u_{E_1}$. Figure~\ref{fig: bif}(a) shows these $E_n$ branches, along with other equilibria that we computed by applying the bifurcation analysis software MatCont \cite{Dooge2003} to the Galerkin truncation \eqref{eq: f} with $N=32$. Each $E_n$ branch has maximum mean energy of $\E\approx3.2067$ when $L\approx1.1947n$. Many of the bounds on $\ov\E$ reported in the next section are approximately equal to the upper envelope of the $E_n$ curves in figure~\ref{fig: bif}(a). In such cases we conclude that the bounds are sharp and are saturated by one of the primary equilibria.

The bifurcation structure of Galerkin truncations of the KSE informs the number of modes that we must include in our bounding computations. With $N$ fixed, the deviation of the ODE system from the full PDE worsens as $L$ is increased. In the PDE, each $E_n$ state is $n$ exact copies of the $E_1$ state. In a truncation with fixed $N$, the branches become more under-resolved as $n$ increases; the spectral content of the $E_n$ branch is every $n^\text{th}$ sine mode, so $E_n$ is approximated by no more than $N/n$ modes. This under-resolution causes the ODE system to be unbounded when $L\ge N/2+1$ because every $E_n$ branch for $n\in [N/2+1,N]$ is truncated after only the single mode $u_n$, and the linear instability of the $u_n$ mode does not saturate without the $u_{2n}$ mode. Thus, for a given $L$, finite bounds are possible for the truncated system or the full KSE only if $N\ge 2L-2$.

\section{\label{sec: results}Bounds on time-averaged energy}

We have computed upper bounds on $\ov\E$ both for large Galerkin truncations of the KSE and for the full KSE by solving the SOS optimization problems \eqref{eq: B trunc kse} and \eqref{eq: B pde kse}, respectively. The software YALMIP \cite{Lofberg2004, Lofberg2009} was used to reformulate the SOS problems as SDPs and interface with SDP solvers. When possible we used the solver MOSEK \cite{Mosek}. As SDPs grew in size because of increasing polynomial degree or number of modes, the memory footprint of MOSEK became prohibitive (e.g., well over 128 GB). For these large SDPs we used the solver SCS \cite{ocpb:16, scs}, whose first-order method converges to a given precision much more slowly than MOSEK's interior-point method but requires less memory. Appendix \ref{app: num} details our computational implementation, including convergence criteria.
\begin{table}[t]
	\caption{\label{tab: conv}Upper bounds on mean energy for odd solutions of the KSE in a domain of size $2\pi L=7\pi$, computed using quadratic and quartic $V$. Values tabulated for various numbers of modes ($N$) are the energy of the third primary equilibrium ($\E_{E_3}$) in the truncated system, bounds on mean energy in the truncated system ($B_{N,2d}$) computed using \eqref{eq: B trunc kse}, and bounds on mean energy in the full PDE ($B^{pde}_{N,2d}$) computed using  \eqref{eq: B pde kse}. To the tabulated precision, $\E_{E_n}=B_{N,4}$ for each $N$. Missing $B_{N,2}^{pde}$ and $B_{N,4}^{pde}$ values did not converge and may be infinite.}
	\centering
	\begin{tabular}{cc l c ll c ll}
		&& &&\multicolumn{2}{c}{quadratic $V$} && \multicolumn{2}{c}{quartic $V$}\\
		\cline{5-6} \cline{8-9}
		$N$ && 
		\multicolumn{1}{c}{$\E_{E_3}$} && 
		\multicolumn{1}{c}{$B_{N,2}$} & 
		\multicolumn{1}{c}{$B^{pde}_{N,2}$} && 
		\multicolumn{1}{c}{$B_{N,4}$} & 
		\multicolumn{1}{c}{$B^{pde}_{N,4}$} \\
		\hline
		6   && 3.126\,656 	    && 9.074\,245 &   			&& 3.126\,656 	&	\\
		7   && \textquotedbl 	&& 9.929\,391 &   			&& \textquotedbl 	& 	 \\
		8   && \textquotedbl 	&& 8.767\,367 &   			&& \textquotedbl 	& 	\\
		9   && 3.173\,413    	&& 9.013\,554 & 34.427\,55 && 3.173\,413	&12.370\,65 \\
		10 && \textquotedbl 	&& 8.986\,130 & 12.015\,03 && \textquotedbl & 3.546\,301\\
		11 && \textquotedbl		&& 9.002\,760 & 10.481\,22 && \textquotedbl & 3.451\,817\\
		12 && 3.174\,051 	    && 9.002\,593 & 9.312\,328 && 3.174\,051    & 3.177\,370\\
		13 && \textquotedbl		&& 9.003\,391 & 9.173\,322 && \textquotedbl & 3.176\,181\\
		14 && \textquotedbl		&& 9.003\,369 & 9.054\,233 && \textquotedbl & 3.176\,139\\
		15 && \textquotedbl		&& 9.003\,411 & 9.030\,412 && \textquotedbl & 3.174\,076\\
		16 && \textquotedbl		&& 9.003\,408 & 9.010\,593 && \textquotedbl & 3.174\,071\\
		17 && \textquotedbl		&& 9.003\,410 & 9.007\,983 && \textquotedbl & 3.174\,068\\
		18 && \textquotedbl		&& 9.003\,409 & 9.004\,389 && \textquotedbl & 3.174\,054\\
		\hline
	\end{tabular}
\end{table}

For an example of how bounds converge as the number of modes $N$ increases, consider the KSE in a domain with $L=3.5$. At this $L$, the $E_3$ branch has the largest mean energy among all equilibria shown in figure~\ref{fig: bif}(a). For each $N$ between 6 and 18, table~\ref{tab: conv} reports $\E_{E_3}$ in the truncated system, as well as bounds for the truncated system and for the full PDE that we computed using both quadratic and quartic auxiliary functionals. The value of $\E_{E_3}$ converges quickly and is accurate to 7 digits in the 12-mode system, wherein the $E_3$ branch is resolved by 4 nonzero modes. In the case of quadratic $V$, the truncated bounds $B_{N,2}$ and PDE bounds $B^{pde}_{N,2}$ both converge toward the optimal bound \eqref{eq: V opt quad} of the background method, which is not a sharp bound on $\ov\E$ at this $L$. In the case of quartic $V$, on the other hand, $B_{N,4}$ and $B^{pde}_{N,4}$ both converge to $\E_{E_3}$ up to 6 digits. This suggests that the $E_3$ state maximizes $\ov\E$ among odd solutions of the KSE at $L=3.5$, and that a quartic $V$ provides the sharp bound.

It is better to compute $B_{N,2d}$, as opposed to $B^{pde}_{N,2d}$, when one's objective is to approximate the large-$N$ limit, thereby approximating the best bound provable within some infinite-dimensional class $\V$ of auxiliary functionals. This is because $B_{N,2d}$ is less expensive to compute at each $N$ (cf.\ Appendix~\ref{app: num}) and often converges faster as $N\to\infty$, as in the quartic-$V$ bounds of table~\ref{tab: conv}. We have computed such results up to the largest domain sizes in which it was tractable to approximate the large-$N$ limit. Section \ref{sec: trunc} reports these findings.

On the other hand, one must compute $B^{pde}_{N,2d}$, as opposed to $B_{N,2d}$, when one's objective is to obtain rigorous bounds for the PDE. The exact value of $B^{pde}_{N,2d}$ at any finite $N$ constitutes such a bound. Full rigor requires also accounting for error in the formulation of the SOS program as an SDP, and in the numerical solution of the SDP itself. This can be done with interval arithmetic using the software VSDP \cite{Jansson2006}, as illustrated in \cite{Goluskin2018}. In the present work we simply approximate $B^{pde}_{N,2d}$ numerically as a proof of concept for the PDE bounding framework of \S\ref{sec: non-quad pde}. Section \ref{sec: pde} reports these findings.

\subsection{\label{sec: trunc}Bounds for ODE truncations}

We have computed bounds $\ov\E\le B_{N,2d}$ for the truncated KSE by solving \eqref{eq: B trunc kse} using truncated auxiliary functionals $V(\ba)$ of polynomial degrees 2, 4, and 6. To approximate the large-$N$ limit at a given value of $L$, we determined the primary equilibrium $E_n$ with the largest mean energy and then used the criterion $N\ge3n$. In computationally easier cases, including all bounds with quadratic $V$, we included many more than $3n$ modes. Results in these over-resolved cases suggest that the $N\ge3n$ criterion gives values of $B_{N,2d}$ within 1\% of the large-$N$ limit, $B_{2d}$, as reflected in the example of table~\ref{tab: conv} where the criterion requires $N\ge9$. Bounds were computed with quadratic $V$ for $L\le100$, with quartic $V$ for $L\le26.6$, and with sextic $V$ for $L\le10$.

\begin{figure}[t]
\centering
\begin{tikzpicture}
\node at (0,0) {\includegraphics[height=140pt,trim={34pt 0 50pt 16pt},clip]{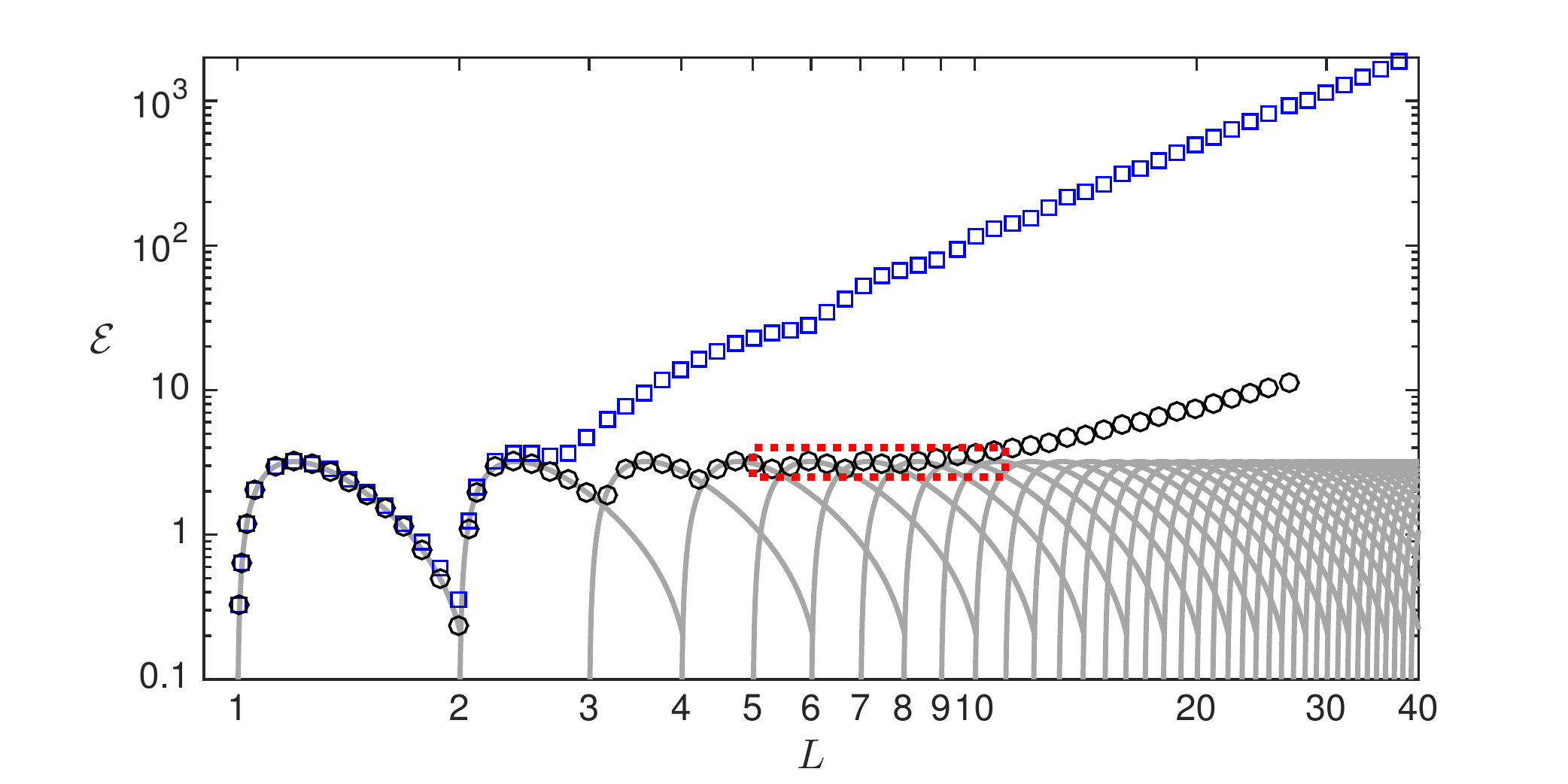}};
\node at (-3.2,1.9) {(a)};
\end{tikzpicture}
\hspace{0pt}
\begin{tikzpicture}
\node at (0,0) {\includegraphics[height=140pt,trim={13pt 0 28pt 16pt},clip]{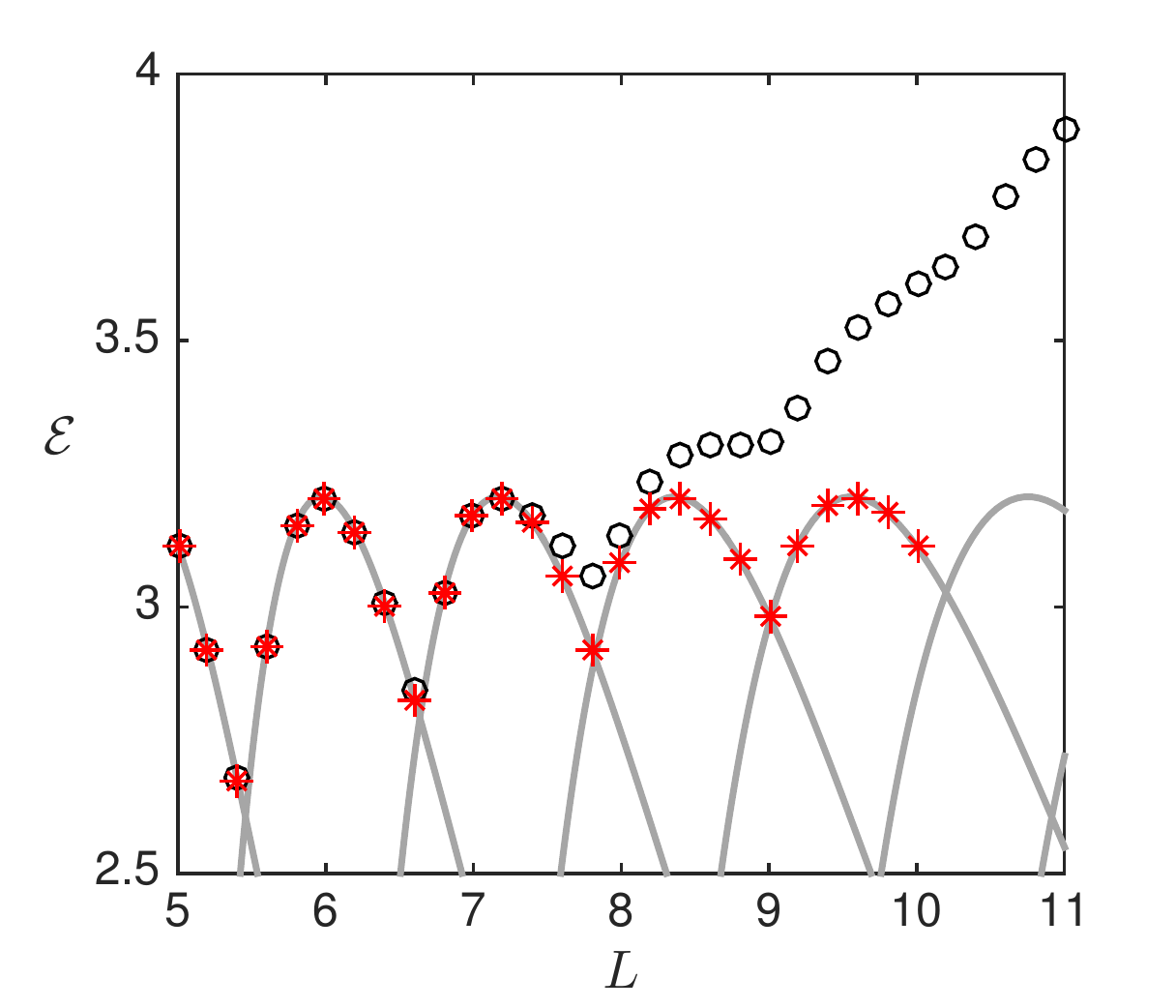}};
\node at (-1.5,1.9) {(b)};
\end{tikzpicture}
\caption{\label{fig: trunc all}
Upper bounds on the mean energy of odd solutions, computed for large truncations of the KSE using auxiliary functionals of degree 2 ({\color{blue}\tiny $\mathbf{\square}$}), 4 ($\circ$), and 6 ({\color{red}$*$}). Mean energies of primary equilibria are shown also ({\color{gray}\solidrule}). Degree-6 bounds are shown only in panel (b), which is a detailed view of the boxed region in panel (a).}
\end{figure}

\begin{figure}[t]
\centering
\includegraphics[height=110pt,trim={95pt 5pt 95pt 10pt},clip]{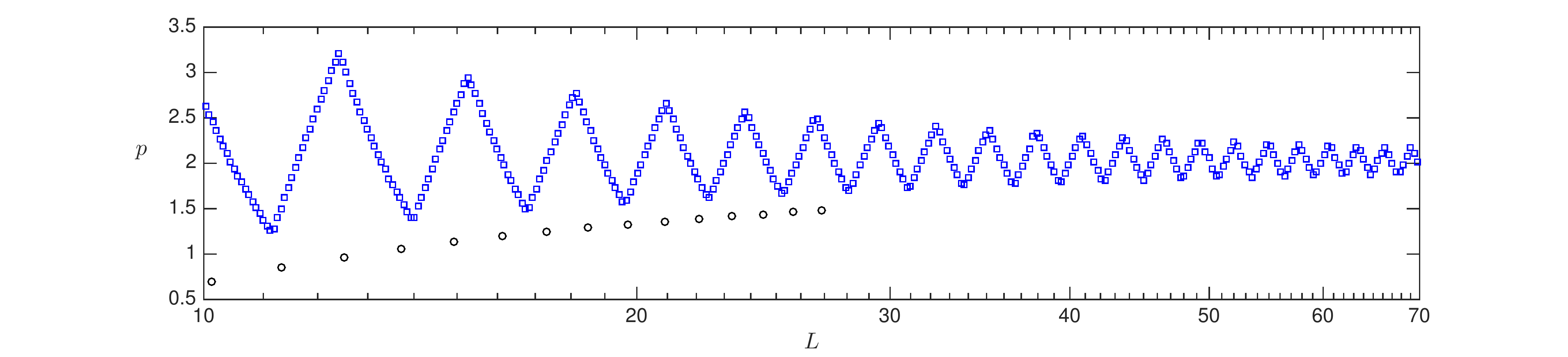}
\caption{\label{fig: exponents}
Local exponents of power laws $c L^p$ fit to upper bounds on mean energy for auxiliary functionals of degree 2 ({\color{blue}\tiny $\mathbf{\square}$}) and 4 ($\circ$). The bound values are shown in figure \ref{fig: trunc all}(a) over a different range of $L$.}
\end{figure}

Bounds on $\ov\E$ for various $L$, computed using quadratic and quartic $V$, are plotted in figure~\ref{fig: trunc all}(a). Energies $\E_{E_n}$ of the KSE's primary equilibria are plotted also for comparison. A region of this plot is expanded in figure~\ref{fig: trunc all}(b), which shows bounds computed using quartic and sextic $V$. When a bound is saturated by one of the $E_n$ branches (up to numerical error) we conclude that it is sharp.

The polynomial degree of $V$ is evidently very important. Quadratic $V$ give sharp bounds (to within 1\%) only for $L\lesssim1.25$, whereas quartic $V$ give sharp bounds for $L\lesssim7.5$. Sextic $V$ give sharp bounds over the entire range $L\le10$ for which we computed them. We suspect that $V$ of any fixed polynomial degree $2d$ can give sharp bounds on $\ov\E$ only up to some finite domain size.

Even if higher-degree polynomial $V$ cannot produce $O(1)$ bounds when $L\gg1$, they might yield bounds that scale better than those proved using quadratic $V$. Our quadratic-$V$ bounds, which are tantamount to the optimal background method, are fit well by $B_2\sim 1.12\,L^2$ over $40\le L\le100$. (Appendix~\ref{app: background} gives examples of the corresponding background functions $\zeta$.) This is the same scaling as the bound $\E_\infty\le O(L^2)$ produced analytically \cite{Bronski2006} and computationally \cite{Fantuzzi2015} using the background method. Apparently bounding $\ov\E$ directly instead of $\E_\infty$ cannot improve this scaling, although our prefactor is more than 20 times smaller. The quartic-$V$ bounds in figure~\ref{fig: trunc all} have not reached the asymptotic regime, but we can discern a trend by examining their local slope.

After refining the results of figure~\ref{fig: trunc all}(a) by computing bounds at many more values of $L$, we have fit local power laws of the form $c\hspace{1pt}L^p$. That is, we approximate the exponent $p$ by a local linear fit of $\log B_{2d}$ to $\log L$. Figure~\ref{fig: exponents} shows how $p$ varies with $L$ for quadratic and quartic $V$. In the quadratic case, the optimal bounds have a local exponent that oscillates while converging to~2. The oscillations have a period of approximately 2.8 with no obvious connection to the KSE's bifurcation structure; they might correspond to bifurcations in the Euler--Lagrange equations governing the optimal background function $\zeta$ in \eqref{eq: V opt quad}. In the case of quartic $V$, on the other hand, bounds oscillate mildly with a shorter period that reflects the shape of the envelope of the $E_n$ energies, as can be seen in figure~\ref{fig: trunc all}(b). This oscillation is not evident in figure~\ref{fig: exponents} because, to make the trend clearer, we have estimated $p$ using quartic-$V$ bounds computed only at values of $L$ where the envelope of $E_n$ energies has a local maximum. That is, we include bounds computed at $L=1.1947\,n$ for $8\le n\le 23$. The limit of $p$ is apparently larger than $2/3$, meaning that when $L\gg1$ the quartic-$V$ bounds will not be as good as the $O(L^{2/3+})$ bounds proved by entropy methods \cite{Otto2009, Goldman2015}. It remains unclear whether the limit of $p$ for quartic-$V$ bounds is less than 2, which would indicate an improvement upon the $O(L^2)$ scaling of quadratic-$V$ bounds.

\subsection{\label{sec: pde}Bounds for the full PDE}

\begin{figure}[t]
\centering
\begin{tikzpicture}
\node at (0,0) {\includegraphics[width=280pt,trim={44pt 36pt 54pt 14pt},clip]{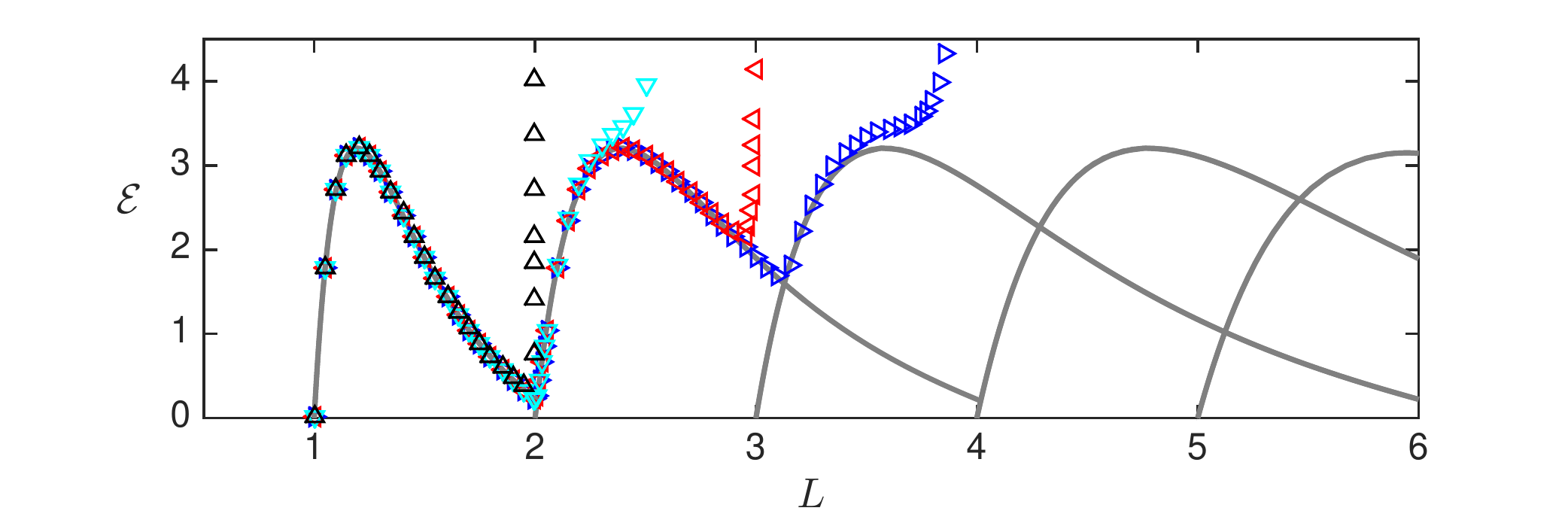}};
\node at (-3.8,1.06) {(a)};
\end{tikzpicture}
\\
\begin{tikzpicture}
\node at (0,0) {\includegraphics[width=280pt,trim={44pt 4pt 54pt 14pt},clip]{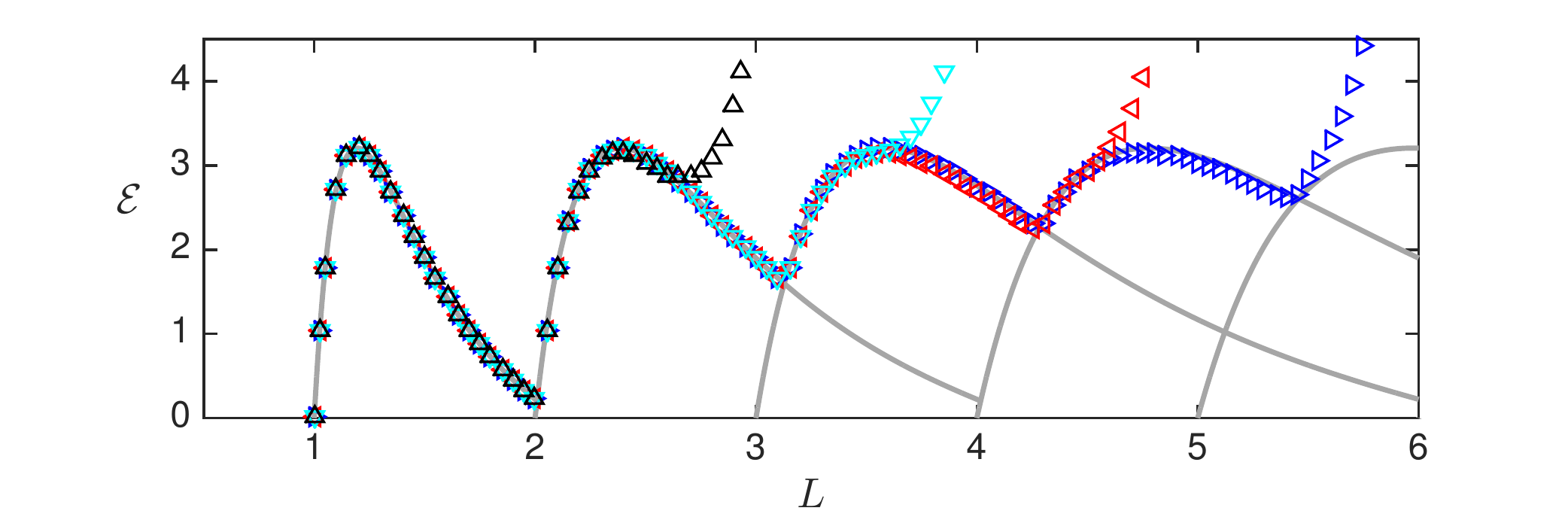}};
\node at (-3.8,1.37) {(b)};
\end{tikzpicture}
\caption{\label{fig: N}Upper bounds on mean energy computed using quartic $V$: (a) bounds $B_{N,4}^{pde}$ for the full KSE and (b) bounds $B_{N,4}$ for truncations of the KSE. The number of modes ($N$) used to compute the bounds are 4 ($\vartriangle$), 6 ({\color{cyan}$\triangledown$}), 8 ({\color{red}\large $\triangleleft$}), and 10 ({\color{blue}\large $\triangleright$}). Mean energies of primary equilibria are shown also ({\color{gray}\solidrule}).}
\end{figure}

We have computed bounds $\ov\E\le B^{pde}_{N,4}$ for the KSE by solving \eqref{eq: B pde kse} for quartic auxiliary functionals $V(\ba,q^2)$. These bounds cannot be sharp when $L\gtrsim7.5$ since they cannot be better than the quartic-$V$ bounds for the truncated system shown in figure~\ref{fig: trunc all}(b), but we find they can be sharp at smaller $L$.

Figure~\ref{fig: N}(a) shows bounds $B^{pde}_{N,4}$ computed for various $L$ using $N=4$, 6, 8, and 10 modes. Every plotted point is a valid bound for the PDE, or at least it would be if the SDP producing it were solved to infinite precision. Bounds are very good for small domains and then quickly blow up to infinity as $L$ increases. When more modes are included, this blowup is postponed until larger $L$. For each $N$, the bounds $B^{pde}_{N,4}$ blow up at smaller $L$ than the corresponding bounds $B_{N,4}$ for the truncated system, which are shown in figure~\ref{fig: N}(b). The latter blow up as $L\to N/2+1$, which is when the $N$-mode truncated system becomes unbounded (cf.\ \S\ref{sec: steady}). Unlike the PDE bounds in figure~\ref{fig: N}(a), some $B_{N,4}$ values plotted in figure~\ref{fig: N}(b) are \emph{not} valid bounds for the PDE---they lie slightly below the $E_n$ curves. This occurs because under-resolved $E_n$ states in truncations of the KSE have slightly less mean energy than in the full KSE. 

\begin{figure}[t]
\centering
\includegraphics[width=280pt,trim={44pt 4pt 54pt 8pt},clip]{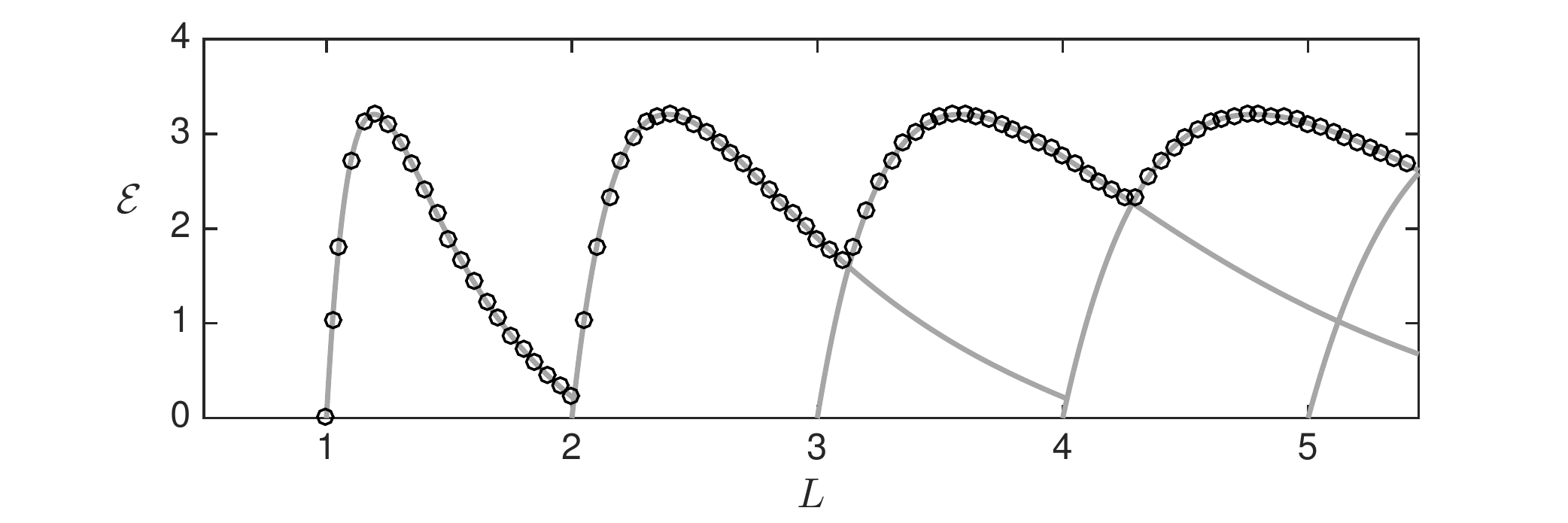}
\caption{\label{fig: pde}Upper bounds on $\ov\E$ for the KSE computed by solving \eqref{eq: B pde kse} with quartic $V$. The bounds saturated by each of the first four branches of $E_n$ equilibria ({\color{gray}\solidrule}) were computed using $N = 6n$ modes.}
\end{figure}

Figure~\ref{fig: pde} shows nearly perfect bounds on $\ov\E$ that we have computed for the KSE by solving \eqref{eq: B pde kse}. The number of modes included in the auxiliary functional \eqref{eq: V pde} to produce bounds saturated by the $E_1$, $E_2$, $E_3$, and $E_4$ branches were 6, 12, 18, and 24, respectively.

\section{\label{sec: con}Conjectures and conclusions}

We have illustrated methods for bounding time averages in nonlinear PDEs by constructing polynomial auxiliary functionals. This is a generalization of the background method, which is tantamount to using a subset of quadratic auxiliary functionals. Polynomial auxiliary functionals and the resulting bounds can be constructed with computer assistance. This is done by formulating sufficient conditions as polynomial optimization problems subject to sum-of-squares constraints, translating these problems into semidefinite programs, and solving the latter computationally. Two related approaches have been presented, both of which involve approximating the PDE by a system of ODEs derived by Galerkin truncation. The first is to compute bounds for the ODE systems using existing methods \cite{Chernyshenko2014, Fantuzzi2016, Goluskin2018} and raise the truncation order until bounds converge. The limit is a bound for the PDE, provided that all PDE solutions are limits of solutions of truncated systems. The second approach is similar but also incorporates analytical estimates on the deviation between the PDE and its truncations. This produces bounds applying to the full PDE, despite the truncations being finite.

We have applied these methods to the Kuramoto--Sivashinsky equation, computing upper bounds on the spatiotemporally averaged energy $\ov\E$ of odd solutions for a variety of domain sizes $2\pi L$. Using the first approach where bounds are computed for large ODE truncations, we have obtained bounds using auxiliary functionals $V$ of degrees 2, 4, and 6. For each domain size, the number of modes was increased until the bounds converged. With $V$ of fixed degree, bounds appear sharp for domains up to a certain size, beyond which they become increasingly conservative as $L$ grows. Quadratic $V$ give sharp bounds (to within 1\%) only for $L\lesssim1.25$, while quartic $V$ give sharp bounds for $L\lesssim7.5$, and sextic $V$ give sharp bounds on the entire range $L\le10$ that was computationally accessible. The reason sharp bounds can be identified as such is that they are saturated (up to numerical error) by one of the primary equilibria, $E_n$, each of which is simply $n$ copies of the first nonzero equilibrium.

To demonstrate the second approach, where bounds for the full PDE are produced by augmenting ODE truncations with analytical estimates, we used quartic $V$ to compute sharp bounds on $\ov\E$ for $L\lesssim5.5$. This approach is more computationally expensive and may not be applicable to as broad a range of nonlinear PDEs, but it is merited when one seeks results that are rigorous to the standard of a computer-assisted proof. This could be accomplished by augmenting SDP computations with interval arithmetic as in \cite{Goluskin2018}.

Our findings bear on the conjecture that solutions of the KSE obey $\ov\E\le O(1)$ for $L\gg1$. The sharp bounds we computed for $L\le10$ suggest the following stronger conjecture.

\begin{conj}
\label{conj: steady}
For all odd solutions $u(x,t)$ of the KSE \eqref{eq: ks} and $L>1$, 
\beq
\max_{u(x,t)}\ov\E = \max_n\E_{E_n}.
\eeq
\end{conj}

\noindent
As $L\to\infty$ the righthand maximum is achieved by periodic copies of the $E_1$ state with the largest mean energy---the equilibrium with spatial period $2\pi L\approx 7.48$ that is shown in figure~\ref{fig: bif}(c). Conjecture \ref{conj: steady} could be strengthened in two ways: by not restricting $u(x,t)$ to the odd subspace, and by replacing the time average $\ov\E$ with $\E_\infty$, the maximum energy at late times. We are not aware of counterexamples to these stronger conjectures; the computational results in figures 14 and 15 of \cite{Cvitanovic2010b} are consistent with them. However, our bounds give positive evidence only for the weaker conjecture stated above, and only for $L\le10$. Regarding the use of auxiliary functionals to prove Conjecture \ref{conj: steady}, our findings suggest another conjecture. 

\begin{conj}
\label{conj: V}
(a) For odd solutions $u(x,t)$ of the KSE \eqref{eq: ks} and any fixed $L$, there exists an auxiliary functional $V$ of finite polynomial degree $2d\in2\mathbb{N}$ that can be used to prove Conjecture 1. In particular, as the truncation order $N$ approaches infinity, bounds on $\ov\E$ for the truncated system given by \eqref{eq: B trunc kse} and bounds for the full PDE given by \eqref{eq: B pde kse} satisfy, respectively,
\begin{align}
\lim_{N\to\infty} B_{N,2d} &= \max_n\E_{E_n},
&
\lim_{N\to\infty} B^{pde}_{N,2d} &= \max_n\E_{E_n}.
\label{eq: conj 2}
\end{align}
(b) However, satisfying \eqref{eq: conj 2} as $L\to\infty$ requires $2d\to\infty$, meaning that $V$ of fixed polynomial degree cannot provide sharp bounds in the large-domain limit.
\end{conj}

The $L$-dependence of our bounds hints at what sort of auxiliary functionals might provide $O(1)$ bounds when $L\gg1$. Apparently the best bounds provable using quadratic $V$ are $\ov\E\le O(L^2)$. Increasing the polynomial degree of $V$ improves bounds on $\ov\E$ at fixed $L$, but it is unclear whether $V$ of any fixed degree can give bounds growing more slowly than $O(L^2)$ when $L\gg1$. As asserted by Conjecture~\ref{conj: V}, it seems that sharp bounds require the degree of $V$ to continue growing as $L\to\infty$. If so, then $O(1)$ upper bounds cannot be proved analytically using any polynomial $V$ ansatz. It remains possible, however, that a relatively simple $V$ ansatz with non-polynomial dependence on $u$ can produce $O(1)$ bounds.

Beyond the KSE, the background method has been used to bound time averages governed by the Navier--Stokes equations and related PDEs. In many cases the resulting bounds are not sharp. Our findings for the KSE suggest that many such results of the background method could be sharpened by generalizing auxiliary functionals beyond the simplest quadratic case. This can be done analytically or computationally. Here we have demonstrated computational approaches that can give sharp bounds, although only when the intrinsic dimension of the dynamics is low enough for the computations to be tractable.

Pushing our SDP-based methods to higher-dimensional regimes requires overcoming the relatively poor scalability of algorithms for solving SDPs. Here we tackled fairly large SDPs by using a solver that implements a first order algorithm, rather than an interior point algorithm. This reduces memory requirements but greatly increases computation time. A second option is to use more tractable relaxations of sum-of-squares constraints based on linear or second-order cone programming~\cite{Ahmadi2017b}, and a third is to use methods specialized for SOS optimization rather than general algorithms for SDPs \cite{Papp2018}. Algorithmic improvements may also be accompanied by restrictions to auxiliary functionals that yield polynomial optimization problems with a favorable sparse structure, to which sparse SOS conditions can be applied~\cite{Waki2006, Lasserre2006}. It remains to be seen how much these techniques can improve upon bounds proved by the background method for a variety of PDEs. For the KSE, at least, our findings demonstrate that the improvement can be quite substantial.

\subsection*{Acknowledgements}
During this work one author (DG) was partially supported by a Van Loo Postdoctoral Fellowship and NSF award DMS--1515161 at the University of Michigan, and by NSERC Discovery Grant RGPIN-2018-04263 and Discovery Accelerator Supplement RGPAS-2018-522657 at the University of Victoria. The other author (GF) was partially supported by EPSRC studentship award 1864077 and by an EPSRC Doctoral Prize Fellowship. Computational resources were provided by Advanced Research Computing at the University of Michigan, and by WestGrid and Compute Canada. Both authors received support from the Geophysical Fluids Dynamics program at the Woods Hole Oceanographic Institution, where some of this work was carried out. We thank Charles Doering, Ian Tobasco, and Andrew Wynn for helpful conversations, and an anonymous referee for insightful comments.

\appendix

\section{\label{app: prop}Symmetry and boundedness of auxiliary functionals}

The following propositions can be used to anticipate properties of auxiliary functionals that produce optimal bounds on $\ov\Phi$ in \eqref{eq: V opt}. Proposition~\ref{prop: sym} asserts that optimal $V$ inherit symmetries from the governing equations. It is formulated for the case of an ODE with a finite group of linear symmetries---a setting that includes the symmetry $a_n\mapsto(-1)^na_n$ of the truncated KSE that is relevant to our present computations. Generalizations to PDEs and Lie symmetry groups are left for future work. Proposition~\ref{prop: bounded} gives conditions under which $V$ must be bounded below. Both propositions are applied to the KSE in \S\ref{sec: poly af kse}.

\begin{prop}
\label{prop: sym}
Consider a well posed ODE $\ddt \ba=\f(\ba)$ with $\ba(t)\in\R^n$. Suppose the ODE and the function $\Phi:\R^n\to\R$ are invariant under a linear transformation $\Lambda:\R^n\to\R^n$, meaning $\f(\Lambda\ba)=\Lambda \f(\ba)$ and $\Phi(\Lambda\ba)=\Phi(\ba)$. Suppose that $\Lambda$ generates a finite symmetry group, meaning $\Lambda^K$ is the identity for some positive integer $K$. Then if there exists any $V:\R^n\to\R$ proving a bound $\ov\Phi\le B$ via the sufficient condition \eqref{eq: S}, there exists such a $V$ that has the symmetry $V(\Lambda\ba)=V(\ba)$.
\end{prop}

\noindent
\emph{Proof.} Suppose $V$ satisfies the sufficient condition \eqref{eq: S}, which in the ODE case is
\beq
S_V(\ba)=B-\Phi(\ba)-\f(\ba)\cdot\nabla V(\ba) \ge 0.
\label{eq: S prop}
\eeq
Consider a symmetrized version of $V$ defined as
\beq
\widehat V(\ba) := \frac1K\sum_{k=0}^{K-1}V(\Lambda^k\ba).
\label{eq: V sym}
\eeq
The invariance $\widehat V(\Lambda\ba)=\widehat V(\ba)$ holds since $\Lambda^K$ is the identity. The claim is proven if $S_{\widehat V}(\ba) \ge 0$ for all $\ba\in\R^n$, and this condition is equivalent to
\beq 
B - \Phi(\ba) - \f(\ba)\cdot \left [\frac1K\sum_{k=0}^{K-1}(\Lambda^k)^T\nabla V(\Lambda^k\ba)\right] \ge0. 
\label{eq: S sym}
\eeq
To show that the above inequality holds, we evaluate \eqref{eq: S prop} at $\Lambda^k\ba$ and use the symmetries of $\f$ and $\Phi$ to find
\beq
B - \Phi(\ba) - \f(\ba)\cdot [(\Lambda^k)^T\nabla V(\Lambda^k\ba)] \ge0.
\eeq
Averaging the lefthand expression for $k=0,1,\ldots,K-1$ implies \eqref{eq: S sym}, thereby proving the claim.

\begin{prop}
\label{prop: bounded}
Suppose the PDE~\eqref{eq: dyn sys} is well posed with solutions $u(x,t)$ remaining in a function space $\U$, meaning that $u(\cdot,t)\in\U$ for all $t\ge0$. Suppose also that $V$ proves a finite bound $\ov\Phi\le B$ via the auxiliary functional condition $B-\Phi(u)-D_V(u)\ge0$, and that $|\Phi(u)|<\infty$ implies $|V(u)|<\infty$. If $\Phi_\infty$ is finite for all PDE solutions $u(x,t)$, then $V(u)$ is bounded below for all $u\in\U$.
\end{prop}

\noindent
\emph{Proof.} If $\Phi_\infty$ defined by \eqref{eq: Phi_inf} is finite there exists an absorbing set $\mathcal A=\{ u\in\U : \Phi(u)\le B_\infty\}$, and we can choose $B_\infty>B$ without loss of generality. Since $\Phi$ is uniformly bounded on $\mathcal A$, by assumption $V$ is also, so $\inf_{\mathcal A}V(u)$ is finite. To show that $V$ is bounded below on $\U$, it suffices to show $\inf_{\U\setminus\mathcal A}V(u)\ge \inf_{\mathcal A}V(u)$.

Let $u_0(x)\in\U\setminus\mathcal A$. Consider a PDE solution $u(x,t)$ with initial condition $u(x,0)=u_0(x)$. There exists a finite first time $T$ at which $\Phi(u(x,T))=B_\infty$, prior to which $\Phi(u(x,t))>B_\infty$. Recalling that $D_V(u(\cdot,t))=\ddt V(u(\cdot,t))$ along all trajectories, we integrate the inequality $\ddt V(u(\cdot,t)) \le B-\Phi(u(\cdot,t))$ up to time $T$ to find
\beq
V(u(x,T)) - V(u_0(x)) \le \int_0^T[B-\Phi(u(x,t))]{\rm d}t \le \int_0^T[B_\infty-\Phi(u(x,t))]{\rm d}t.
\eeq
The righthand integrand is negative, so $V(u_0(x))> V(u(x,T))\ge \inf_{\mathcal A}V(u)$. This proves the claim that $V(u)$ is bound below uniformly for all $u\in\U$. $\square$

\section{\label{app: num}Computational implementation}

When using YALMIP to reformulate the SOS optimization problems \eqref{eq: B trunc kse} and \eqref{eq: B pde kse} as SDPs, we exploit the fact that optimal $V$ are invariant under $a_n\mapsto(-1)^na_n$, as explained in \S\ref{sec: poly af kse}. The polynomial $S$ shares this symmetry, so the Gram matrix representation of $S$, which is central to the translation of SOS constraints into SDP constraints, can be decomposed into two blocks \cite{Gatermann2004, Lofberg2009}. Ordinarily YALMIP can detect symmetries of $S$ and block diagonalize its Gram matrix accordingly, but many of the problems solved here were too large for this functionality. Thus we modified YALMIP to impose the particular block diagonal structure arising in the present cases.

All SDP solutions were computed on a single core. The SDPs were solved by an interior-point algorithm using MOSEK  when the memory footprint to do so was not prohibitively large (e.g., larger than 128 GB). All bounds for the full KSE reported in \S\ref{sec: pde} were computed using MOSEK. Bounds for truncations of the KSE reported in \S\ref{sec: trunc} were computed using MOSEK for all quadratic $V$, for quartic $V$ when $N\le 36$, and for sextic $V$ when $N\le16$. Solver tolerances were set to $10^{-12}$, resulting in many computations terminating when progress stalled. In the worst such cases, relative infeasibilities were still smaller than $10^{-7}$. These infeasibilities alone do not give an error bound on the approximate optimum of the SDPs---that is, on $B_{N,2d}$ or $B^{pde}_{N,2d}$. This would require more sophisticated analysis, as implemented by the software VSDP \cite{Jansson2006}. However, the fact that various bounds are sharp to 5 or more digits (e.g., in table \ref{tab: conv}) suggests that values of $B_{N,2d}$ and $B^{pde}_{N,2d}$ computed using MOSEK are at least this precise. 

The first-order solver SCS was applied to SDPs too large for MOSEK. Solver tolerances were set to $10^{-4}$, and all solutions reached these tolerances. For the largest SDPs solved here, where $V$ was quartic and $N$ was larger than 50, meeting these tolerances required running SCS for over a month at each value of $L$. The formulation of these large SDPs using YALMIP took between a few hours and a day. Comparing SCS and MOSEK solutions for smaller SDPs suggests that tolerances of $10^{-4}$ in SCS give approximate values of $B_{N,2d}$ that are precise to at least 4 digits. This can be seen in table~\ref{tab: solvers}, which shows examples of bounds computed with both solvers. The table also shows the number of iterations required for solutions to meet tolerances, the time of the SDP solver running with a single thread on a 2.5 GHz Intel Xeon processor, and the memory required to formulate the SDP using YALMIP and then run the solver. For the purposes of this table and the next one, all MOSEK tolerances were $10^{-8}$. Although times of different solvers cannot be directly compared, even when tolerances are nominally identical, it is clear that SCS converges more slowly than MOSEK, and that its memory footprint grows more slowly with the problem size.

\begin{table}[t]
\caption{\label{tab: solvers}Upper bounds on mean energy for odd solutions of the truncated KSE, computed using quartic auxiliary functionals ($2d=4$). Tabulated values give the number of modes in the truncation ($N$), the domain size ($2\pi L$), the computed bound ($B_{N,4}$), the SDP solver used, the number of iterations, the time of the SDP solution, and the memory required to set up and solve the SDP.}
\centering
\begin{tabular}{cccclll}
$N$ & $L$ & $B_{N,4}$ & solver & steps & time & memory  \\
\hline
12 & 3 &		1.900\,280 & MOSEK & 17 & 1 sec.  & 0.42 GB \\
&&			1.900\,284 & SCS & 197\,260 & 311 & 0.41 \\
20 & 5 &		3.113\,014 & MOSEK & 17 & 59 & 0.85 \\
&&			3.113\,034 & SCS & 265\,880 & 2\,842 & 0.43 \\
28 & 7 & 		3.174\,205 & MOSEK & 32 & 2\,659 & 5.12 \\
&&			3.174\,393 & SCS & 7\,851\,400 & 304\,200 & 0.63
\end{tabular}
\end{table}

For a fixed number of modes, computing bounds for the truncated KSE by solving \eqref{eq: B trunc kse} is less expensive that computing bounds for the full KSE by solving \eqref{eq: B pde kse}. Table~\ref{tab: trunc vs pde} gives some examples of the wall time required to solve the SDPs arising from \eqref{eq: B trunc kse} and \eqref{eq: B pde kse}, respectively, using MOSEK with a single thread.

\begin{table}[t]
\caption{\label{tab: trunc vs pde}Time required by MOSEK to solve the SDP formulations of \eqref{eq: B trunc kse} and \eqref{eq: B pde kse} with quartic auxiliary functionals to obtain bounds for the truncated KSE ($B_{N,4}$) and full KSE ($B^{pde}_{N,4}$), respectively, for various numbers of modes ($N$) and domain sizes ($2\pi L$).}
\centering
\begin{tabular}{ccclccl}
&& \multicolumn{2}{c}{Solving \eqref{eq: B trunc kse}}
	&& \multicolumn{2}{c}{Solving \eqref{eq: B pde kse}}\\
 \cline{3-4} \cline{6-7}
$N$ & $L$ & $B_{N,4}$ & time && $B^{pde}_{N,4}$ & time \\
\hline
8	& 2	& 0.222\,904	& 0.1	 sec.	&& 0.222\,938 & 0.6 sec.  \\
12 	& 3 	& 1.900\,280 	& 1.0		&& 1.900\,440 & 5.5  \\
16	& 4	& 2.762\,934	& 13.3	&& 2.763\,049 & 41.9
\end{tabular}
\end{table}

\section{\label{app: background}Optimal background functionals}

The bounds we have produced using quadratic auxiliary functionals can be translated into the language of the background method if desired. When $V(\ba)$ is quadratic, the optimization \eqref{eq: B trunc kse} of bounds for the truncated KSE is a truncated version of the background method optimization \eqref{eq: V opt quad} in spectral space. So long as enough modes are included in the truncation, optimal quadratic $V(\ba)$ provide close approximations to the optimal coefficients $\alpha$ and background functions $\zeta(x)$ solving \eqref{eq: V opt quad}. The optimal $V$ have only even-index $a_n$ in their linear terms, so they take the form
\begin{equation}
V(\ba) = c|\ba|^2 + \sum_{n=1}^{\lfloor N/2\rfloor} c_{2n}a_{2n},
\label{eq: V spectral 1}
\end{equation}
and the SDP solver returns optimal values of the coefficients $c$ and $c_{2n}$. Meanwhile, after dropping the irrelevant term $\tfrac{\alpha}{2}\fint{\zeta^2 {\rm d}x}$, the auxiliary functional \eqref{eq: V background} defined using the background function can be expanded in terms of Fourier coefficients as
\begin{equation}
V(u) = \frac{\alpha}{4\pi L}|\ba|^2 - \frac{\alpha}{2\pi L}\sum_{n=1}^\infty a_nz_n,
\label{eq: V spectral 2}
\end{equation}
where
\begin{equation}
\zeta(x) = (\pi L)^{-1/2}\sum_{n=1}^\infty z_n\sin(nx/L).
\end{equation}
When $u(x,t)$ and $\zeta(x)$ in the background method are truncated after $N$ modes, expressions \eqref{eq: V spectral 1} and \eqref{eq: V spectral 2} can be equated to find
\begin{align}
\alpha &= 4\pi L c, & 
\zeta(x) &= -\frac{1}{2 c (\pi L)^{1/2}}
\sum_{n=1}^{\lfloor N/2\rfloor} c_{2n}\sin(2nx/L).
\label{eq: recover bg}
\end{align}

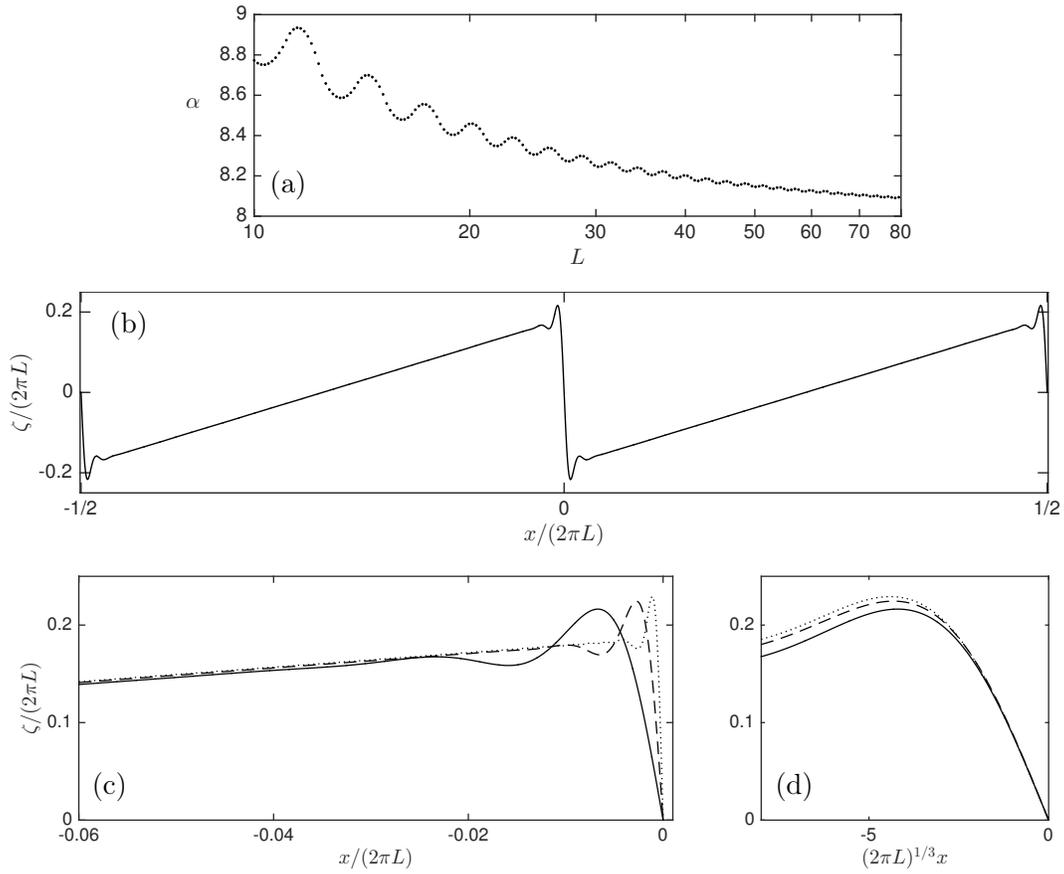
\begin{figure}[tp]
\centering
\hspace{6pt}
\begin{tikzpicture}
\node at (0,0) {\includegraphics[height=100pt,trim={28pt 5pt 50pt 5pt},clip]{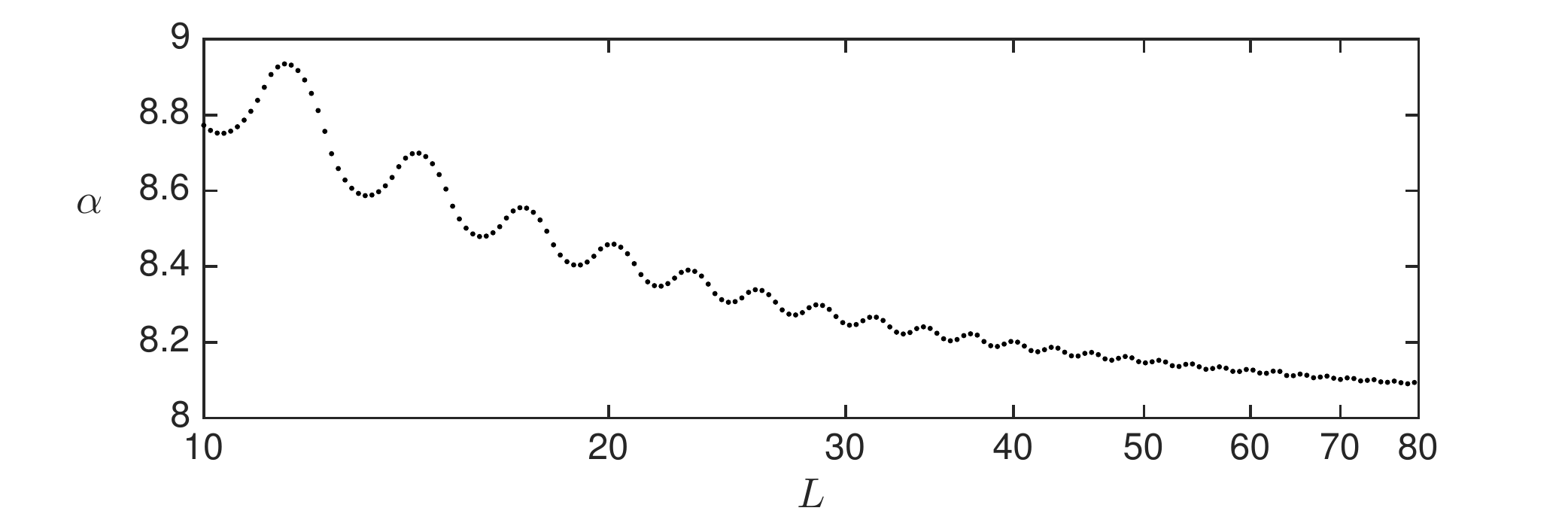}};
\node at (-3.45,-.7) {(a)};
\end{tikzpicture}
\\[2pt]
\begin{tikzpicture}
\node at (0,0) {\includegraphics[height=96pt,trim={60 2 70 15},clip]{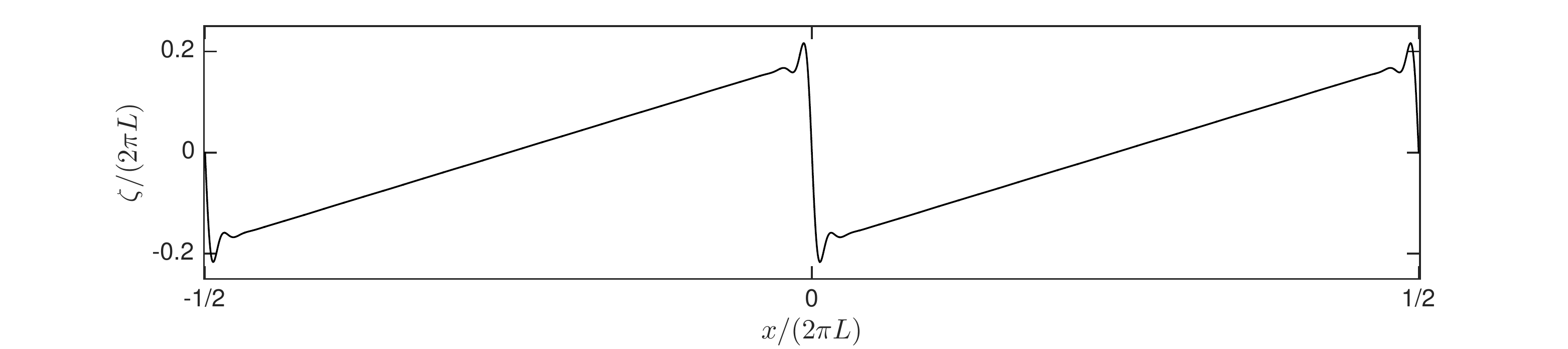}};
\node at (-5.4,1.25) {(b)};
\end{tikzpicture}
\\[0pt]
\begin{tikzpicture}
\node at (0,0) {\includegraphics[height=115pt,trim={25 0 50 12},clip]{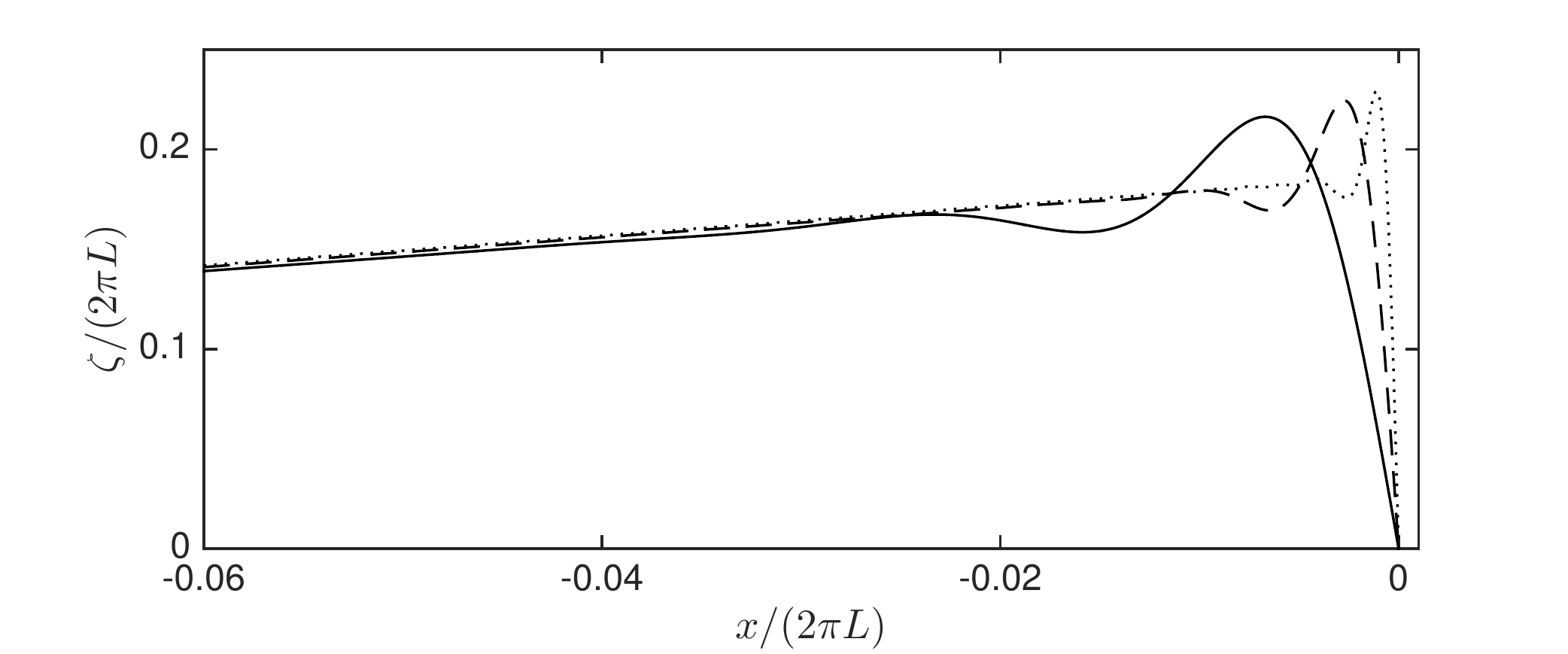}};
\node at (-3.15,-.9) {(c)};
\end{tikzpicture}
\hspace{2pt}
\begin{tikzpicture}
\node at (0,0) {\includegraphics[height=115pt,trim={19 0 20 12},clip]{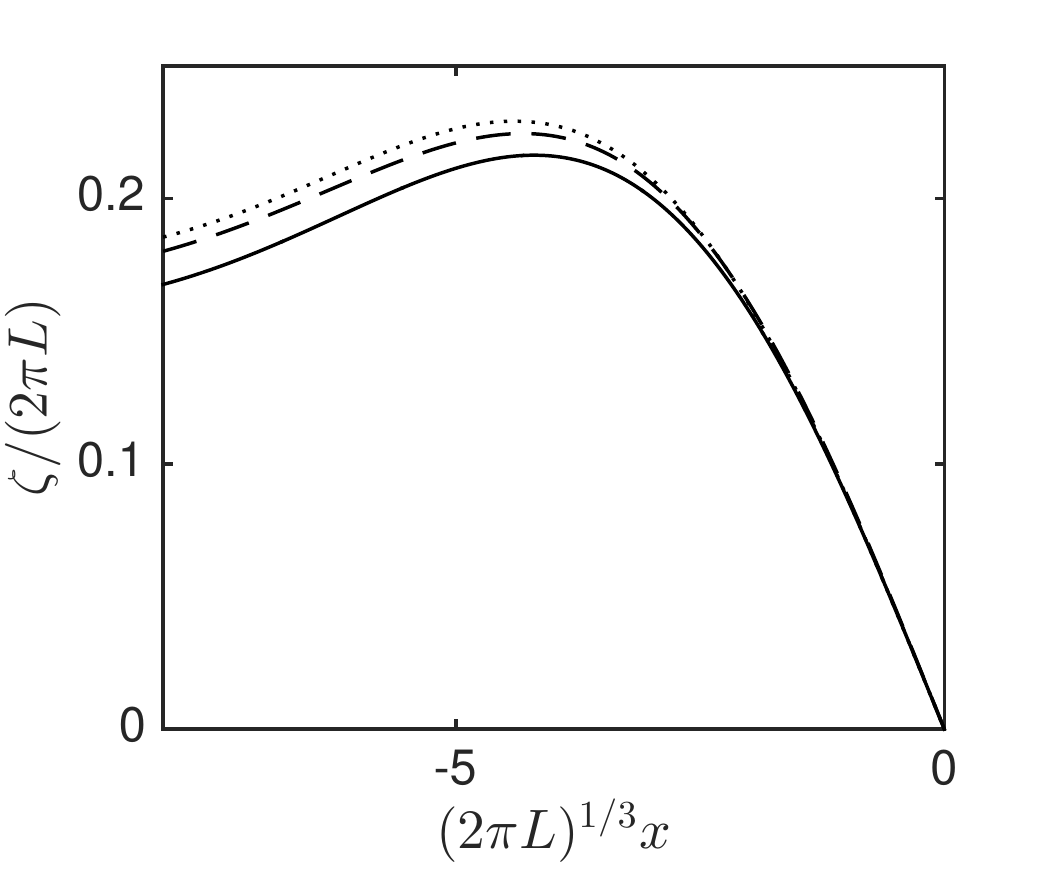}};
\node at (-1.3,-.9) {(d)};
\end{tikzpicture}
\caption{\label{fig: background}
Optimal coefficients $\alpha$ and functions $\zeta(x)$ for the background method \eqref{eq: V opt quad} with various $L$, computed by solving the truncated optimization \eqref{eq: B trunc kse} with $2d=2$ and $N\ge8L$ modes. Displayed quantities are (a) optimal $\alpha$ for various domain sizes $2\pi L$, (b) optimal $\zeta(x)$ for $L=10^{1.3}$, and (c-d) details of optimal $\zeta(x)$ for $L=10^{1.3}$ (\solidrule), $10^{1.6}$ (\longdashedrule), and $10^{1.9}$ (\dottedrule).}
\end{figure}

Figure~\ref{fig: background} shows some optimal $\alpha$ and $\zeta(x)$ for the background method, recovered from our quadratic $V$ according to \eqref{eq: recover bg}. Panel (a) shows that, as $L$ increases, the optimal leading coefficient $\alpha$ oscillates and decreases, seeming to approach a constant value. The oscillations in $\alpha$ coincide with the oscillations in the slope of upper bounds shown in figure~\ref{fig: exponents}, although only the former quantity appears to depend smoothly on $L$. Figure~\ref{fig: background}(b) shows the optimal background function for $L=10^{1.3}\approx19.95$, with both $x$ and $\zeta$ normalized by the domain size $2\pi L$. The shape of $\zeta$ strongly resembles background functions that have been constructed in \cite{Bronski2006, Fantuzzi2015} for the related task of bounding the instantaneous energy $\E_\infty$. Because the optimal $V$ is symmetric under $u(x,t)\mapsto-u(x+\pi L,t)$ (cf.\ \S\ref{sec: poly af kse}), the optimal $\zeta$ has a fundamental period of $\pi L$, although $\zeta$ with a fundamental period of $2\pi L$ also can give $O(L^2)$ bounds \cite{Bronski2006}. The scaling of our optimal $\zeta$ as $L\to\infty$ is similar to that in \cite{Bronski2006, Fantuzzi2015}. Outside the boundary layers the slope $\zeta'(x)$ approaches a constant. This can be seen in figure~\ref{fig: background}(c), which shows normalized $\zeta$ near the central boundary layer for three different domain sizes. The limiting slope appears to be 3/4, as opposed to 9/4 in the $\E_\infty$ bounding problem \cite{Fantuzzi2015}. The boundary layer thickness scales as $O(L^{-1/3})$. This can be seen by the collapse of profiles in figure~\ref{fig: background}(d), where the profiles in panel (c) are re-plotted with $x$ normalized according to this scaling.

\bibliographystyle{plain}
\bibliography{KS_bounds_v12.bbl}

\end{document}